\documentclass[11pt,a4paper]{article}

\usepackage[margin=1in]{geometry}
\usepackage{amsmath,amssymb,amsthm,mathtools,mathrsfs}
\usepackage{xcolor,microtype,indentfirst}
\usepackage{needspace}
\usepackage[numbers,sort&compress]{natbib}
\usepackage[colorlinks=true,linkcolor=blue,citecolor=blue,urlcolor=blue]{hyperref}
\usepackage{newtxtext,newtxmath,cleveref}

\usepackage{indentfirst}
\usepackage{tikz}
\mathtoolsset{showonlyrefs}

\newtheorem{theorem}{Theorem}[section]
\newtheorem{lemma}[theorem]{Lemma}

\newtheorem{remark}{Remark}
\numberwithin{equation}{section}
\counterwithin*{parentequation}{section} 
\newcommand{\R}{\mathbb R}
\newcommand{\D}{\mathbb D}
\newcommand{\Hh}{\mathcal H}
\newcommand{\pz}{\partial_z}
\newcommand{\pb}{\partial_{\bar z}}
\newcommand{\ip}[2]{\langle #1,#2\rangle}
\DeclareMathOperator{\diver}{div}
\DeclareMathOperator{\curl}{curl}

\allowdisplaybreaks[2]

\title{\bfseries\Large Strict hot spots of second Neumann eigenfunctions on simply connected Lipschitz domains in two-dimensional space forms\footnote{\footnotesize The work is supported by the National Natural Science Foundation of China (No.12001276, No.12326303, No.12431018).}}
\author{\normalsize Haiyun Deng$^{1}$, Xuyong Jiang$^{2}$, Xiaoping Yang$^{3}$, Jun Zou$^{3}$\footnote{\footnotesize E-mail: hydeng@nau.edu.cn(H. Deng), jiangxy@cczu.edu.cn(X. Jiang), xpyang@nju.edu.cn(X. Yang), jun\_zou@smail.nju.edu.cn(J. Zou)}\\[4pt]
\small $^{1}$Department of Applied Mathematics, Nanjing Audit University, Nanjing, 211815, China\\
\small $^{2}$Department of Mathematics, Changzhou University, Changzhou, 213164, China\\
\small $^{3}$Department of Mathematics, Nanjing University, Nanjing, 210093, China}
\date{}
\begin{document}
	\maketitle
	
	\begin{abstract}
    We prove that every second Neumann eigenfunction has no interior
    critical points on any bounded simply connected Lipschitz domain
    $\Omega$ in either $\mathbb R^2$ or $\mathbb H^2$.
    The same conclusion holds for such a domain
    $\Omega\subset\mathbb S^2$ if its closure lies in an open hemisphere
    and $\mu_2(\Omega)>2$.
    Consequently, in each case the maximum and minimum are attained
    only on the boundary, thereby verifying the \textit{hot spots
    conjecture} for these three classes of domains.
    The main difficulty is that a Neumann eigenfunction on a
    Lipschitz domain need not belong to $H^2(\Omega)$, so we represent
    its rotated gradient in disk coordinates as the sum of an $H^1$ part
    and an antiholomorphic part with finite weighted $L^2$ norm.
    The key estimate, a strict two-multiplier inequality on the unit
    disk, shows that an interior critical point would contradict the
    nonnegativity of the divergence--curl energy of the comparison fields.
	\end{abstract}

	

	\setcounter{tocdepth}{2}
	\tableofcontents
	
	\section{Introduction}\label{intro}

Let $M_K$ be the simply connected complete two-dimensional space form
of curvature $K\in\{0,-1,1\}$, where $M_0=\R^2$,
$M_{-1}=\mathbb H^2$, and $M_1=\mathbb S^2$ is the unit sphere.
For a bounded Lipschitz domain $\Omega\subset M_K$, we consider
the following Neumann eigenvalue problem
\begin{equation}\label{eq:neumann}
 \begin{cases}
 -\Delta u=\mu u & \text{in }\Omega,\\
 \partial_n u=0 & \text{on }\partial\Omega.
 \end{cases}
\end{equation}
All geometric operators refer to the
standard metric on $M_K$, and $n$ is the outward unit normal, defined
almost everywhere on $\partial\Omega$.
Since the boundary need not be smooth, the equation is
understood in the weak sense: $u\in H^1(\Omega)$ satisfies
\begin{equation}\label{weak}
	\int_\Omega\nabla u\cdot\nabla\varphi d A
	=\mu\int_\Omega u\varphi d A,
	\qquad \varphi\in H^1(\Omega).
\end{equation}
For a bounded Lipschitz domain $\Omega$, the Rellich compactness theorem
gives the compact embedding $H^1(\Omega)\hookrightarrow L^2(\Omega)$,
and hence the self-adjoint Dirichlet and Neumann Laplacians have discrete
spectra. 
The eigenvalues, counted with multiplicity, are denoted by
\begin{equation}\label{spectrum}
	0=\mu_1(\Omega)<\mu_2(\Omega)\le\mu_3(\Omega)\le\cdots.
\end{equation}
Throughout, $\mu=\mu_2(\Omega)$ denotes the second Neumann eigenvalue, and we write $\lambda_D=\lambda_1^D(\Omega)$ for the first
Dirichlet eigenvalue.

Rauch’s hot spots conjecture \cite{Rauch1974} asserts that any eigenfunction of the first positive Neumann eigenvalue of a planar domain attains its maximum and minimum only on the boundary; precise formulations are given in \cite{BanuelosBurdzy1999JFA}.  For general domains in $\mathbb{R}^n$, the conjecture is false.  For instance, even for planar domains with interior holes, Neumann eigenfunctions can attain their maxima inside the region. The conjecture has been proven false for certain planar multiply connected domains \cite{Burdzy1999Ann,Burdzy2005Duke}, a bounded starlike Lipschitz domain in $\R^3$ \cite{AtarBurdzy2026Arxiv} and convex domains in sufficiently high dimensions \cite{deDios2024Arxiv}. For planar domains, the hot spots conjecture has been proved for several	classes, including convex domains with two axes of symmetry \cite{JerisonNadirashvili2000JAMS}, lip domains \cite{Atar2024JAMS}, nearly circular convex domains \cite{Miyamoto2009JMP}, and triangles \cite{Judge2020Ann,Judge2022Ann}. Further results concern critical points and monotonicity on triangles
\cite{ChenGuiYao}, and symmetry, multiplicity, and critical points on isosceles trapezoids and kites \cite{DengEtAlKites}. Methods used in these works include probabilistic coupling \cite{BanuelosBurdzy1999JFA}, domain deformation, and nodal set analysis. 

Recently, OpenAI \cite{OpenAI2026} stated a theorem excluding interior critical points of every second Neumann eigenfunction on smooth bounded simply connected planar domains, without a simplicity assumption. Taking the result of OpenAI as an input, Burdzy \cite{BurdzyBoundary2026} derived the non-strict hot-spots property on bounded simply connected planar domains with compact Neumann resolvent and simple first positive Neumann eigenvalue, without assuming boundary smoothness; and Ba\~nuelos \cite{Banuelos2026Arxiv} extended the hot spots conclusion to simply connected domains of finite area with compact Sobolev embedding and simple first positive Neumann eigenvalue.  Independently of \cite{OpenAI2026}, Deng, Jiang and Yang \cite{DengJiangYang20262} developed a near-boundary rigidity theorem and proved that  the maximum modulus of every second Neumann eigenfunction on a bounded convex planar domain is attained only on the boundary, without assuming boundary smoothness or eigenvalue simplicity. Moreover, they also proved strict boundary maximum modulus for higher Neumann eigenfunctions under suitable conditions, including all third Neumann eigenfunctions
on arbitrary triangles. Pineda Montoya et al. \cite{Pineda2026Arxiv} proved that every second eigenfunction on piecewise smooth Jordan domains  with convex corners has no interior critical points.  For quantitative comparisons between interior and boundary extrema, see \cite{DeDiosPontHsuTaylor,DengJiangYang20261,Mariano2023PA,Steinerberger2023RMI}. For spaces with constant curvature, Hatcher \cite{Hatcher2026Arxiv} successfully adapted Miyamoto's method in \cite{Miyamoto2009JMP} to convex domains in the hyperbolic plane $\mathbb{H}^2$. Hatcher showed that if the second Neumann eigenvalue $\mu_2(\Omega) \leq 1/4$ for a convex domain $\Omega\subset\mathbb{H}^2$, then the second Neumann eigenfunction has no interior critical points. For related results, see \cite{Chen2019JMPA,Chen2026AMPA,Deng2025AMPA,Hatcher2025SIAMJMA,Judge2025AMQ,Judge2022CPDE,Kennedy2026TAMS,Langford2023MathAnn} and references therein.

Our main result excludes interior critical points for every second Neumann eigenfunction, without assuming that $\mu$ is simple.
	
\begin{theorem}\label{main}
Let $\Omega\subset M_K$ be a bounded simply
 connected Lipschitz domain, where $K\in\{0,-1,1\}$.
 If $K=1$, assume in addition that
 $\overline\Omega$ is contained in an open hemisphere and
 $\mu_2(\Omega)>2$. Then every second Neumann eigenfunction $u$ satisfies
 \begin{equation}\label{noncritical}
  \nabla u(p)\ne0\quad\text{ for every }p\in\Omega.
 \end{equation}
 Consequently, the following strict hot spots results hold:
 \begin{equation}\label{extrema}
 	\min_{\partial\Omega}u<u(p)<\max_{\partial\Omega}u
 	\quad\text{ for every }p\in\Omega.
 \end{equation}
\end{theorem}

{\noindent \bf Idea of the proof:}
Suppose that a second Neumann eigenfunction $u$ has a critical point
$p_0$, and let $W=J\nabla u$, where $J$ is rotation by $\pi/2$.
On a Lipschitz domain $u\in H^1(\Omega)$ may fail to belong to
$H^2(\Omega)$, so its rotated gradient is a priori only in
$L^2(\Omega;T\Omega)$. Its divergence and curl nevertheless satisfy
$\diver W=0$ and $\curl W=-\mu u$ in distributions, and the Neumann
condition gives $W$ zero weak tangential trace. Thus $W$ belongs to
the energy space $\mathscr X$. The Hodge decomposition and the
spectral gap $\mu<\lambda_D$ give $E_\mu(Y)\ge0$ for every
$Y\in\mathscr X$, with $E_\mu(W)=0$ (\Cref{hodge}). We use the zero
of $W$ at $p_0$ to construct two comparison fields whose combined
divergence--curl energy would be strictly smaller than $E_\mu(W)$.

We choose a conformal map $\Psi:\D\to\Omega$ with $\Psi(0)=p_0$ and
write $\Psi^*dA=\rho d z$. The rescaled  disk coefficient $F_W$ vanishes
at the origin, and \eqref{gradnorms} gives
$\sqrt\rho\,F_W,\pz F_W\in L^2(\D)$ without providing a global
bound for $\pb F_W$. The Cauchy transform produces an $H^1$ function
with the same $\pz$ derivative as $F_W$, leaving an antiholomorphic
remainder. In the plane this gives the required decomposition of $U=F_W$.
On the curved surfaces we use a Hilbert-space-valued lift to account
for the factor $\kappa$ inside $\pb(\kappa\pz F_W)$. In all three geometries we obtain a solution
of $-\Delta U=VU$ with $U(0)=0$ and
\[
 U=U_{\rm a}+\widehat U,\qquad
 \pz U_{\rm a}=0,\qquad
 \sqrt\rho\,U_{\rm a}\in L^2(\D;\Hh),\qquad
 \widehat U\in H^1(\D;\Hh).
\]
Thus the derivative term of $\mathcal A_V$ involves only $\widehat U$.
For the potential term, the pointwise bound
$\|V^{1/2}U_{\rm a}\|^2\le C\rho\|U_{\rm a}\|^2$ shows why weighted
$L^2$ control of $U_{\rm a}$ suffices.

The key step is the strict two-multiplier inequality of \Cref{strict}.
We multiply the trace of $\widehat U$ by $\gamma_1=\cos\theta$ and
$\gamma_2=\sin\theta$ to prescribe the regular boundary data for two
solutions $U_1,U_2$ of the same disk equation. Their antiholomorphic
parts are obtained by multiplying $U_{\rm a}$ by the corresponding
extensions of $\gamma_i$, whose poles at $0$ are canceled by
$U_{\rm a}(0)=0$. To compare the three energies, write $U=h+q_h$,
where $h$ is harmonic and $q_h\in H^1_0(\D;\Hh)$, and decompose $U_i$
in the same way. Multiplication by $\gamma_1\pm i\gamma_2=e^{\pm i\theta}$
shifts the Fourier indices by one, so the two harmonic comparison
functions contribute $2\pi\|h(0)\|^2$ more than $h$ to the $\pz$ term
of the disk energy, by \eqref{harm-energy}. Since $U(0)=0$, the Green
formula for $-\Delta q_h=VU$ gives
$h(0)=-q_h(0)=-\int_\D G(0,t)V(t)U(t) d t$.
With this value of $h(0)$, completing squares expresses $-Q_V$ as
a sum of nonnegative terms in \eqref{exactsquares}. Strict positivity
of the kernel in \Cref{poskernel} then gives
\[
 Q_V:=\mathcal A_V(U_1)+\mathcal A_V(U_2)-\mathcal A_V(U)<0
 \qquad\text{if }V^{1/2}U\ne0 .
\]
In the Euclidean and hyperbolic planes,
$\|V^{1/2}U\|_{L^2(\D;\Hh)}^2=\mu\|\nabla u\|_{L^2(\Omega)}^2>0$.
On the sphere, curvature changes the factor $\mu$ to $\mu-2$, which
explains the hypothesis $\mu>2$.

To return to $\Omega$, we use smooth multipliers $\beta_i$ with
boundary values $\gamma_i$ such that $U_i-\beta_iU\in H^1_0(\D;\Hh)$.
This relation transfers the weak tangential boundary condition of $W$
to both comparison fields. In the Euclidean plane their energies
agree with the disk energies. In the hyperbolic plane, projection
preserves the potential term and can only decrease the derivative
term. On the sphere, projection leaves additional boundary terms,
which cancel in $Q_V$ because $\gamma_1^2+\gamma_2^2=1$ on the
circle. Thus in each geometry $Q_V$ is at least the sum of the two
nonnegative divergence--curl energies, contradicting the strict disk
inequality.

\vspace{1em}
{\noindent\bf Comparison with previous work:}
Similar to OpenAI \cite{OpenAI2026}, we use the nonnegativity of the
divergence--curl energy and boundary multipliers on the disk.
 OpenAI constructs a family of multipliers from a
reciprocal Green kernel and shows that the existence of an interior
critical point would force each of the corresponding boundary energies
to vanish. We use the two fixed multipliers $\cos\theta$ and
$\sin\theta$. For this pair, \Cref{strict} shows that an interior
critical point would force the sum of the energies of the two
comparison solutions to be strictly less than the energy of the
original disk solution.
Unlike \cite{Banuelos2026Arxiv} and \cite{BurdzyBoundary2026}, they take OpenAI's result as an external input and obtain the corresponding results through smooth approximation. Our proof works directly on $\Omega$, independently of OpenAI's smooth-domain theorem, and does not use any approximation by smooth domains.

\vspace{1em}
The rest of this paper is organized as follows.  In Section \ref{spectral} we use the strict eigenvalue comparison
$\mu<\lambda_D$ to prove that the divergence--curl energy $E_\mu$ is
nonnegative on $\mathscr X$. Its expression in conformal coordinates
motivates the problem studied in Section \ref{disksec}, where we
prove the strict two-multiplier inequality for Hilbert-valued
solutions of $-\Delta U=VU$ on $\D$. In Section \ref{proofmain} we
apply this inequality to $F_W$ in the Euclidean
plane and to the corresponding lifts in the hyperbolic plane and on
the unit sphere, completing the proof of \Cref{main}.

	\section{Preliminaries and common lemmas}\label{spectral}
	
	Throughout this section, let $M_K$ denote the two-dimensional space form with its standard metric,
	where $K=0,1,-1$ corresponds to the Euclidean plane, the unit sphere,
	and the hyperbolic plane, respectively.
	We say that $\Omega\subset M_K$ is a Lipschitz domain if, for every boundary point $p\in\partial\Omega$, there exists a coordinate chart
	$\Upsilon:U\to B_1(0)$ centered at $p$ and a Lipschitz function
	$\phi:\mathbb R\to\mathbb R$ with $\phi(0)=0$ such that
	$\Upsilon(U\cap\Omega)=B_1(0)\cap D$, where $D$ is the region above the graph of $\phi$.
	Let $J$ denote the orientation-compatible parallel rotation by $\pi/2$ on each tangent plane, and let $n$ and $\tau=Jn$ denote the unit outward normal and unit tangent vector on $\partial\Omega$, respectively, wherever they exist.
	The Lipschitz assumption guarantees that unit outward normal and unit tangent vectors exist almost everywhere with respect to arclength measure.

	For a bounded Lipschitz domain $\Omega$, the standard Rellich compactness theorem gives the compact embedding
	$H^1(\Omega)\hookrightarrow L^2(\Omega)$. Consequently, the self-adjoint Dirichlet and Neumann Laplacians
	have compact resolvents. Each has discrete spectrum and an orthonormal eigenbasis in $L^2(\Omega)$.
	The eigenvalues, counted with multiplicity, are characterized by the min--max principle.
	Let $u$ be a weak Dirichlet or Neumann eigenfunction on a bounded Lipschitz
	domain. Then $u\in H^1_0(\Omega)$ in the Dirichlet
	case and $u\in H^1(\Omega)$ in the Neumann case. In both cases, standard
	interior elliptic regularity gives $u\in C^\infty(\Omega)$.
    Global
	$H^2(\Omega)$ regularity can fail for either boundary condition, for example
	at a reentrant corner.
	The boundary H\"older estimates of Gilbarg and Trudinger \cite[Theorem 8.29]{GT1983}
	for the Dirichlet problem and of Nittka \cite[Theorem 3.14]{Nittka2011JDE}
	for the Neumann problem show that $u$ has global H\"older regularity, that is,
	\begin{equation}\label{regularity}
		u\in C^{0,\alpha}(\overline\Omega)\qquad\text{for some }0<\alpha<1.
	\end{equation}

	Next, we compare $\mu_2(\Omega)$ and $\lambda_1^D(\Omega)$ for bounded Lipschitz domains $\Omega\subset M_K$.
	In the Euclidean plane case, the strict inequality $\mu_2(\Omega)<\lambda_1^D(\Omega)$ is obtained in \cite{Filonov2005}.
	In the spherical case, if $\Omega$ is contained in a hemisphere but is not itself a hemisphere, then
	the strict inequality $\mu_2(\Omega)<\lambda_1^D(\Omega)$ is obtained in \cite{Ashbaugh1997JEDP}.
	In the hyperbolic plane case , the strict inequality $\mu_2(\Omega)<\lambda_1^D(\Omega)$ holds
	for smooth bounded simply connected domains, as established in \cite{Hua2026MA}.
	A simple argument shows that the same strict comparison holds for any bounded Lipschitz domain $\Omega$.
	Indeed, let $B$ denote the geodesic disk with the same area as $\Omega$.
	The hyperbolic Szeg\H{o}--Weinberger inequality \cite[Theorem 6.1]{Ashbaugh1995JLMS}
	and the Faber--Krahn inequality (see \cite{Benguria2007Duke}) yield the inequalities
	\[
	\mu_2(\Omega)\le\mu_2(B)<\lambda_1^D(B)\le\lambda_1^D(\Omega),
	\]
	where the strict middle inequality follows from applying the smooth case to $B$.

	\begin{lemma}\label{lem:spectral-compare}
		Let $\Omega\subset M_K$ be a bounded Lipschitz domain.
		In the spherical case assume that $\Omega$ is compactly contained in
		the hemisphere. Then
		\begin{equation}\label{spectral-gap}
			0<\mu_2(\Omega)<\lambda_1^D(\Omega).
		\end{equation}
	\end{lemma}

	\subsection{The divergence--curl energy}\label{intrinsic}
	
	Let $\Omega\subset M_K$ be a bounded Lipschitz domain, and let
	$L^2(\Omega;T\Omega)$ denote the space of square-integrable tangent vector
	fields on $\Omega$. 
	For $Y\in L^2(\Omega;T\Omega)$, we denote its divergence and curl in the weak sense by $\diver Y$ and $\curl Y$, respectively. More precisely, for every test function $\varphi\in C_0^\infty(\Omega)$, we define
	\begin{equation}\label{dist-div-curl}
		\begin{aligned}
			\langle\diver Y,\varphi\rangle
			=-\int_\Omega\ip{Y}{\nabla\varphi}dA,\qquad
			\langle\curl Y,\varphi \rangle
			=-\int_\Omega\ip{Y}{J\nabla\varphi}dA .
		\end{aligned}
	\end{equation}
	By Cauchy--Schwarz, both pairings are bounded by
	$\|Y\|_{L^2(\Omega)}\|\nabla\varphi\|_{L^2(\Omega)}$, so they extend
	continuously to $H^1_0(\Omega)$, thus they are belong to $H^{-1}(\Omega)$.
	If $Y$ is a smooth vector field,
	these definitions agree with the usual divergence
	and curl.
	For $h\in H^1(\Omega)$, then $Y=J\nabla h\in L^2(\Omega;T\Omega)$. Since $J$ preserves inner products, we have
	\begin{equation}\label{curl-rot-grad}
		\langle\curl(J\nabla h),\varphi \rangle=\langle\diver(\nabla h),\varphi \rangle
		=-\int_\Omega\langle\nabla h,\nabla\varphi\rangle\,dA
		=\langle\Delta h,\varphi\rangle .
	\end{equation}
	Thus $\curl(J\nabla h)=\diver(\nabla h)=\Delta h$ in the weak sense.
	Similarly, 
		$\curl(\nabla h)=\diver(J\nabla h)=0$ in the weak sense.
        
	We write $\curl Y\in L^2(\Omega)$ when there is an $f\in L^2(\Omega)$
	such that $\displaystyle\langle\curl Y,\varphi \rangle=\int_\Omega f\varphi dA$ for every
	$\varphi\in C_0^\infty(\Omega)$, in which case we write $\operatorname{curl} Y = f$. Define
	\begin{equation}\label{tang}
		\mathcal T(Y,\varphi)=
		\int_\Omega(\curl Y)\varphi dA
		+\int_\Omega\ip{Y}{J\nabla\varphi} dA
		\quad \text{ for }\varphi\in H^1(\Omega).
	\end{equation}
	If $Y$ is a smooth vector field, Stokes' formula gives
	\begin{equation}\label{tang-smooth}
		\mathcal T(Y,\varphi)=\int_{\partial\Omega}(Y\cdot\tau)\varphi \mathrm ds  .
	\end{equation}
	For vector fields with square-integrable divergence and curl, 
	in view of \eqref{tang-smooth}, 
	we impose the boundary condition $Y\cdot\tau=0$ weakly through the space
	\begin{equation}\label{fields}
		\mathscr X=\{Y\in L^2(\Omega;T\Omega):
		\diver Y,\curl Y\in L^2(\Omega),\
		\mathcal T(Y,\varphi)=0\text{ for every }\varphi\in H^1(\Omega)\}.
	\end{equation}
	For a vector field with square integrable divergence and curl, we define the divergence–curl energy by
	\begin{equation}\label{intrinsicenergy}
		E_\mu(Y)=\int_\Omega
		\bigl(|\diver Y|^2+|\curl Y|^2-\mu|Y|^2\bigr) dA .
	\end{equation}
	The following lemma establishes the nonnegativity of $E_\mu$ and characterizes its nullspace.
	
	\begin{lemma}\label{hodge}
		Let $\Omega $ be a bounded simply connected Lipschitz domain with $\mu_2<\lambda_D$, and set $\mu=\mu_2$.
		Then for every $Y\in\mathscr X$, we have
		\begin{equation}\label{hodge-bound}
			E_\mu(Y)=\|\diver Y\|_{L^2(\Omega)}^2+\|\curl Y\|_{L^2(\Omega)}^2-\mu\|Y\|_{L^2(\Omega)}^2\ge0 .
		\end{equation}
		For every second Neumann eigenfunction $u$, the vector field $W=J\nabla u\in\mathscr X$, and $E_\mu(W)=0$.  
		Conversely, every field in
		$\mathscr X$ with zero energy is of this form.
	\end{lemma}
	
	\begin{proof}
		We solve the Dirichlet and Neumann Poisson problems in their weak forms,
		\begin{equation}\label{poisson-system}
			\Delta v=\diver Y,\quad v=0 \text{ on } \partial\Omega,\qquad\qquad
			\Delta h=\curl Y,\quad\int_\Omega h\,dA=0,\quad\partial_n h=0 \text{ on } \partial\Omega.
		\end{equation}
		Taking $\varphi=1$ in the weak boundary condition \eqref{fields} gives $ \int_\Omega\curl Y\,dA=0$,
		which is the compatibility condition for
		the Neumann problem.
		Hence, there exist $v\in H^1_0(\Omega)$ and $h\in H^1(\Omega)$, with
		$\int_\Omega h\,dA=0$, satisfying
		\begin{equation}\label{poisson-weak}
			\begin{aligned}
				\int_\Omega\ip{\nabla v}{\nabla\varphi} dA
				&=-\int_\Omega(\diver Y)\varphi dA &&(\varphi\in H^1_0(\Omega)),\\
				\int_\Omega\ip{\nabla h}{\nabla\psi} dA
				&=-\int_\Omega(\curl Y)\psi dA &&(\psi\in H^1(\Omega)).
			\end{aligned}
		\end{equation}
		For any $\varphi\in H^1(\Omega)$, using the definition \eqref{tang},
		together with  $\curl \nabla v=0$, $\diver(J\nabla\varphi)=0$ and $\curl(J\nabla h)=\Delta h$ in the weak sense, we have
		\begin{equation}\label{poisson-tangent}
			\begin{aligned}
				\mathcal T(\nabla v,\varphi)
				&=\int_\Omega\ip{\nabla v}{J\nabla\varphi} dA=0,\\
				\mathcal T(J\nabla h,\varphi)
				&=\int_\Omega\bigl((\Delta h)\varphi+
				\ip{\nabla h}{\nabla\varphi}\bigr) dA=0 .
			\end{aligned}
		\end{equation}
		Thus $\nabla v, J\nabla h\in \mathscr{X}$ .

		Set $Z=Y-\nabla v-J\nabla h$. Since both gradient fields belong to
		$\mathscr X$, thus $Z\in \mathscr X$,
		by using
		$\diver(J\nabla h)=\curl(\nabla v)=0$ in the weak sense, we obtain the following equations in the weak sense:
		\begin{equation}\label{harmonic-remainder}
			\diver Z=\diver Y-\Delta v=0,\qquad
			\curl Z=\curl Y-\Delta h=0 \quad\text{in }\Omega .
		\end{equation}
		The one-form $Z^\flat$ is closed and coclosed, so interior elliptic regularity for the Hodge Laplacian yields that $Z$ is smooth in $\Omega$. Since $\diver Z=0$, the one-form $(-JZ)^\flat$ is
		closed, and simple connectedness gives a smooth function $k$ with
		$\nabla k=-JZ$, or equivalently $Z=J\nabla k$. To verify that
		$k\in H^1(\Omega)$, fix an interior ball $B\Subset\Omega$ and truncate
		$k$ to $k_N=\max(-N,\min(k,N))$. The Poincar\'e inequality with the
		mean taken on $B$ yields
		\begin{equation}\label{primitive-bound}
			\|k_N\|_{L^2(\Omega)}
			\le C\bigl(\|\nabla k_N\|_{L^2(\Omega)}+\|k_N\|_{L^2(B)}\bigr)
			\le C\bigl(\|Z\|_{L^2(\Omega)}+\|k\|_{L^2(B)}\bigr).
		\end{equation}
		Since $|\nabla k|=|Z|$, the bound is independent of $N$, and monotone
		convergence gives $k\in L^2(\Omega)$ and hence $k\in H^1(\Omega)$.
		We may therefore test the weak boundary condition with $k$, combain with $\curl Z=0$ and $J\nabla k=Z$, obtaining
		$\mathcal T(Z,k)=\|Z\|_{L^2(\Omega)}^2=0$.
		Thus $Z=0$, and
		\begin{equation}\label{hodge-split}
			Y=J\nabla h+\nabla v,\qquad
			\int_\Omega\ip{\nabla v}{J\nabla h} dA=0 .
		\end{equation}
		Here the decomposition is $L^2$-orthogonal. Indeed, $v\in H^1_0(\Omega)$ and
		$\diver(J\nabla h)=0$ imply
		\begin{align}
			\int_\Omega\langle\nabla v,J\nabla h\rangle dA=0.
		\end{align}

		The Dirichlet and Neumann Laplacians have compact resolvents. Choose
		real $L^2$-orthonormal eigenbases $(\psi_j)_{j\ge1}$ and
		$(\phi_k)_{k\ge1}$ for the Neumann and Dirichlet problems, respectively,
		and write $h_j=\int_\Omega h\psi_j dA$ and
		$v_k=\int_\Omega v\phi_k dA$ to be the coefficients of $h$ and $v$ in
		Neumann and Dirichlet eigenbases.
		Testing \eqref{poisson-weak} with the
		corresponding eigenfunctions gives
		\begin{equation}\label{poisson-coefficients}
			\int_\Omega(\curl Y)\psi_j dA=-\mu_jh_j,\qquad
			\int_\Omega(\diver Y)\phi_k dA=-\lambda_k^D v_k .
		\end{equation}
		Using the orthogonal decomposition \eqref{hodge-split} and testing
		the weak equations \eqref{poisson-weak} with $v$ and $h$, respectively, 
		together with \eqref{poisson-coefficients} and Parseval's identity, we have
		\begin{equation}\label{parseval}
			\begin{aligned}
				\|Y\|_{L^2(\Omega)}^2
				&=\|\nabla v\|_{L^2(\Omega)}^2+\|\nabla h\|_{L^2(\Omega)}^2
				=-\int_\Omega v\,\diver Y dA-\int_\Omega h\,\curl Y dA\\
				&=\sum_{j\ge2}\mu_j|h_j|^2+
				\sum_{k\ge1}\lambda_k^D|v_k|^2,\\
				\|\curl Y\|_{L^2(\Omega)}^2+\|\diver Y\|_{L^2(\Omega)}^2
				&=\sum_{j\ge2}\left|\int_\Omega(\curl Y)\psi_j dA\right|^2
				+\sum_{k\ge1}\left|\int_\Omega(\diver Y)\phi_k dA\right|^2\\
				&=\sum_{j\ge2}\mu_j^2|h_j|^2+
				\sum_{k\ge1}(\lambda_k^D)^2|v_k|^2 .
			\end{aligned}
		\end{equation}
		Here the sum over Neumann eigenfunctions starts at $j=2$ because $h$ has zero mean and
		$\int_\Omega\curl Y dA=0$. Substituting these expansions into
		the definition of $E_\mu$, we obtain
		\begin{equation}\label{exactHodge}
			E_\mu(Y)=\sum_{j\ge2}\mu_j(\mu_j-\mu)|h_j|^2+
			\sum_{k\ge1}\lambda_k^D(\lambda_k^D-\mu)|v_k|^2 .
		\end{equation}
		The factors $\mu_j-\mu$ are nonnegative for $j\ge2$, while
		$\lambda_k^D-\mu>0$ by assumption, so every term in \eqref{exactHodge}
		is nonnegative, which proves \eqref{hodge-bound}.

		For every second Neumann eigenfunction $u$, and $W=J\nabla u$, we have
		$\diver W=0$ and $\curl W=\Delta u=-\mu u$, thus
		\begin{equation}\label{rotated-energy}
			\begin{aligned}
				\mathcal T(W,\varphi)
				&=-\mu\int_\Omega u\varphi dA
				+\int_\Omega\ip{\nabla u}{\nabla\varphi} dA=0,\\
				E_\mu(W)&=\mu^2\|u\|_2^2-\mu\|\nabla u\|_2^2=0 .
			\end{aligned}
		\end{equation}
		Conversely, suppose $Y\in\mathscr X$ and $E_\mu(Y)=0$. Each factor
		$\lambda_k^D(\lambda_k^D-\mu)$ in the Dirichlet sum in \eqref{exactHodge}
		is strictly positive, so $v_k=0$ for every $k$.
		The Neumann sum similarly gives $h_j=0$ for $\mu_j>\mu$, so the
		normalization $h_1=0$ implies $h$ is a second Neumann eigenfunction. The orthogonal
		decomposition \eqref{hodge-split} now gives $Y=J\nabla h$, which proves the characterization
		of the nullspace of $E_\mu$.
	\end{proof}
	
	\subsection{Conformal transport}\label{conformal}

	For $\Omega\subset M_K$ as in \Cref{main}, choose a conformal
	coordinate $\zeta=\xi+i\upsilon$.
	Precisely, we use Cartesian coordinates
	in the plane, Poincar\'e disk coordinates in the hyperbolic plane, and
	stereographic projection from the pole opposite the center of the
	containing hemisphere on the sphere.
	In each case
	$\Omega_{\rm coord}=\zeta(\Omega)$ is a bounded simply connected
	Lipschitz domain. Moreover, in the spherical case and hyperbolic case, the coordinate image is compactly contained
	in the unit disk $\mathbb D$, that is $|\zeta|\le R <1$ for some $R$.
	With this conformal coordinate $\zeta$, we have
	\begin{equation}\label{metric-model}
		g=e^{2\delta_0(\zeta)}|d\zeta|^2, \quad e^{\delta_0(\zeta)}=\frac{1+K^2}{1+K|\zeta|^2}.
	\end{equation}
	
	Given any fix $p_0\in\Omega$,
	since $\Omega_{\rm coord}$ is a bounded simply connected domain, the Riemann mapping theorem guarantees the existence of a Riemann map $\Phi:\mathbb D\to\Omega_{\rm coord}$,
	such that $\Psi=\zeta^{-1}\circ\Phi$ satisfies $\Psi(0)=p_0$.
	We set
	\begin{equation}\label{conformaldata}
		\delta=\delta_0\circ\Phi,\qquad \rho=e^{2\delta}|\Phi'|^2,\qquad
		\kappa=e^{-2\delta}.
	\end{equation}
	Thus $\Psi^*g=\rho|\mathrm dz|^2$, and we use $ d z$ for Euclidean area on $\D$.
	
	A bounded Lipschitz Jordan domain is a quasidisk. Quasiconformal
	reflection \cite{Ahlfors1963Acta} extends its Riemann map $\Phi$ through the circle, and higher
	integrability \cite[Theorem 1]{Gehring1973Acta} on a neighborhood of the closed
	disk gives $\Phi'\in L^{2r}(\mathbb D)$ for some $r>1$.
	The metric coefficient and its first derivatives are bounded on
	$\overline{\Omega_{\rm coord}}$, thus we have
	\begin{equation}\label{metricbounds}
		\Phi'\in L^{2r}(\D),\quad \rho\in L^r(\D),\quad
		e^{\pm\delta}\in W^{1,2r}(\D)\cap L^\infty(\D),\quad
		0<c\le e^\delta\le C .
	\end{equation}
	The exponent and constants may depend on the domain.

	We show that the conformal map $\Psi=\zeta^{-1}\circ\Phi$ induces
	an isomorphism between $H^1(\Omega)$ and $H^1(\mathbb D)$ by pullback.
	For a function $f$ on $\Omega$, write $v=f\circ\Psi$, so conformal
	invariance of the Dirichlet integral and change of variables give
	\begin{align}\label{transport-integrals}
		\int_{\mathbb D}|\nabla v|^2 dz=\int_\Omega|\nabla f|^2 dA,
		\qquad
		\int_{\mathbb D}\rho|v|^2 dz=\int_\Omega|f|^2 dA.
	\end{align}
	If $f\in H^1(\Omega)$, we first show that $\|v\|_{H^1(\mathbb D)}\le C\|f\|_{H^1(\Omega)}$.
	Set $D_0=B(0,1/2)$.
	Since $\Phi$ is smooth on the interior of $\mathbb D$,
	$\rho$ admits a positive lower bound on $\overline{D_0}$,
	so the weighted integral controls $\|v\|_{L^2(D_0)}$. Let $T_Nv$ be its truncation
	to $[-N,N]$. The Poincar\'e inequality with the mean taken on $D_0$
	gives
	\begin{align}\label{pullback-bound}
		\|T_Nv\|_{L^2(\mathbb D)}
		\le C\bigl(\|\nabla T_Nv\|_{L^2(\mathbb D)}
		+\|T_Nv\|_{L^2(D_0)}\bigr)
		\le C\|f\|_{H^1(\Omega)}.
	\end{align}
	Letting $N\to\infty$, we have $\|v\|_{H^1(\mathbb D)}\le C\|f\|_{H^1(\Omega)}$.
	Conversely, for $v\in H^1(\D)$, H\"older's inequality and the
	Sobolev embedding on the disk give
	\begin{align}\label{weightspace}
		\int_{\mathbb D}\rho|v|^2 dz
		\le \|\rho\|_{L^r(\mathbb D)}\|v\|_{L^{2r/(r-1)}(\mathbb D)}^2
		\le C\|v\|_{H^1(\mathbb D)}^2.
	\end{align}
	Together with conformal invariance of the Dirichlet integral, we have $\|f\|_{H^1(\Omega)}\le C\|v\|_{H^1(\D)}$.
	Therefore, $\Psi$ induces an isomorphism between $H^1(\Omega)$ and $H^1(\mathbb D)$ via pullback,
	and the same correspondence also holds for $H^1_0$.
	Moreover, applying these estimates to the coefficients of vector fields,
	we see that the correspondence also holds for vector fields.

	\subsection{Energy in disk coordinates}
	
	For a tangent field $Y$, we next give a formula for $E_\mu$ in the disk coordinate $z=x+iy$. The oriented orthonormal frame $e_1=e^{-\delta_0}\partial_\xi$, $e_2=e^{-\delta_0}\partial_\upsilon$ identifies each tangent plane with $\mathbb C$ via $ae_1+be_2\mapsto a+ib$. For $Y=Y^1e_1+Y^2e_2$, we denote its complex coefficient on
	$\Omega_{\rm coord}$ by $w_0=(Y^1+iY^2)\circ\zeta^{-1}$ and its
	pullback to the disk by $w=w_0\circ\Phi$. Under this identification of tangent vectors with complex numbers in the disk coordinate, we have the following results.
	
	\begin{lemma}\label{diskenergy}
		Let $Y\in\mathscr X$ have disk coefficient $w$, and set $F_Y=e^\delta w$. Then
		\begin{equation}\label{coefenergy}
			E_\mu(Y)=4\int_\D\kappa|\pz F_Y|^2 d z
			-\mu\int_\D\rho\kappa|F_Y|^2 d z,
		\end{equation}
		where $\delta$, $\rho$ and $\kappa$ are given in \eqref{conformaldata}.
	\end{lemma}

	\begin{proof} 
		Write $e^{\delta_0}w_0=P+iQ$ in the coordinate
		$\zeta=\xi+i\upsilon$. Since the metric dual of $Y$ is
		$P\,\mathrm d\xi+Q\,\mathrm d\upsilon$, the 
		divergence and curl are give as
		\[
		(\diver Y+i\curl Y)\circ\zeta^{-1}
		=e^{-2\delta_0}\bigl(P_\xi+Q_\upsilon
		+i(Q_\xi-P_\upsilon)\bigr)
		=2e^{-2\delta_0}\partial_\zeta(P+iQ).
		\]
		Now $F_Y=(P+iQ)\circ\Phi$, thus the chain rule and the definition
		$\kappa=e^{-2\delta}$  and $w=e^{-\delta}F_Y$ gives
		\begin{equation}\label{field-local-identities}
			\begin{aligned}
				& (\diver Y+i\curl Y)\circ\Psi =\frac{2e^{-2\delta}}{\Phi'}\pz F_Y  =\frac{2\kappa}{\Phi'}\pz F_Y,\\
				& |Y|^2\circ\Psi =|w|^2=e^{-2\delta}|F_Y|^2=\kappa|F_Y|^2.
			\end{aligned}
		\end{equation}
		Here the derivatives in \eqref{field-local-identities}
		are understood in the weak sense.
		Since $\Psi^* dA=\rho d z$ and
		$\rho\kappa=|\Phi'|^2$, changing variables in
		\eqref{field-local-identities} gives
		\begin{equation}\label{field-local-identities-div-curl}
			\begin{aligned}
				\|\diver Y\|_{L^2(\Omega)}^2+\|\curl Y\|_{L^2(\Omega)}^2
				&=\int_\D\rho\frac{4\kappa^2}{|\Phi'|^2}
				|\pz F_Y|^2 d z
				=4\int_\D\kappa|\pz F_Y|^2 d z,\\
				\|Y\|_{L^2(\Omega)}^2&=\int_\D\rho\kappa|F_Y|^2 d z.
			\end{aligned}
		\end{equation}
		These formulas and the definition of $E_\mu$ give \eqref{coefenergy}.
	\end{proof}

		The following lemma rewrites the Neumann eigenfunction equation
		in terms of the disk coefficient $F_W$ of the rotated gradient
		$W=J\nabla u$, and expresses the $L^2$ mass and Dirichlet energy
		of $u$ in disk coordinates.
		
		\begin{lemma}\label{lem:rotated-gradient}
			Let $u$ be a second Neumann eigenfunction and $W=J\nabla u\in\mathscr{X}$ have disk
			coefficient $w$. Set $F_W=e^\delta w$. Then
			\begin{equation}\label{rotated-second}
				-4\pb(\kappa\pz F_W)=\mu\rho\kappa F_W
				\qquad\text{in }\D.
			\end{equation}
			where $\delta$, $\rho$ and $\kappa$ are given in \eqref{conformaldata}.
			Moreover, the $L^2$ mass and the Dirichlet energy of $u$ satisfy
			the following identities in disk coordinates,
			\begin{subequations}\label{gradnorms}
				\begin{align}
					\mu^2\|u\|_{L^2(\Omega)}^2
					&=4\int_\D\kappa|\pz F_W|^2 d z,
					\label{rotated-derivative}\\
					\|\nabla u\|_{L^2(\Omega)}^2
					&=\int_\D\rho\kappa|F_W|^2 d z.
					\label{rotated-mass}
				\end{align}
			\end{subequations}
		\end{lemma}
		
		\begin{proof}
			Set $u_0=u\circ\zeta^{-1}$ and $\widetilde u=u\circ\Psi$.
			In the oriented orthonormal frame associated with $\zeta$, the
			complex coefficient of $\nabla u$ is
			$2e^{-\delta_0}\partial_{\bar\zeta}u_0$.
			Since $J$ acts as multiplication by $i$, pulling back the coefficient
			of $W=J\nabla u$ to the disk gives
			\[
			w=2ie^{-\delta}(\partial_{\bar\zeta}u_0)\circ\Phi.
			\]
			Multiplying by $e^\delta$ and applying the chain rule to
			$\widetilde u=u_0\circ\Phi$, we obtain
			\begin{equation}\label{rotated-coefficient}
				F_W=2i(\partial_{\bar\zeta}u_0)\circ\Phi,\qquad
				2\pb\widetilde u
				=2\overline{\Phi'}(\partial_{\bar\zeta}u_0)\circ\Phi.
			\end{equation}
			For the metric $g=e^{2\delta_0}|\mathrm d\zeta|^2$, the Laplacian
			has the coordinate expression
			\[
			\Delta=e^{-2\delta_0}
			(\partial_\xi^2+\partial_\upsilon^2)
			=4e^{-2\delta_0}\partial_\zeta\partial_{\bar\zeta}.
			\]
			Thus the eigenfunction equation in the interior becomes
			$4\partial_\zeta\partial_{\bar\zeta}u_0
			=-\mu e^{2\delta_0}u_0$, so differentiating the expression for
			$F_W$ and using $\kappa e^{2\delta}=1$, we obtain
			\[
			\kappa\pz F_W
			=2i\kappa\Phi'
			(\partial_\zeta\partial_{\bar\zeta}u_0)\circ\Phi
			=-\frac{i\mu}{2}\Phi'\widetilde u.
			\]
			Differentiating this identity in $\bar z$ and using
			\eqref{rotated-coefficient}, together with $\rho\kappa=|\Phi'|^2$,
			yields
			\[
			-4\pb(\kappa\pz F_W)
			=2i\mu\Phi'\pb\widetilde u
			=\mu|\Phi'|^2F_W
			=\mu\rho\kappa F_W,
			\]
			which proves \eqref{rotated-second}. 
			Since $\diver W=0$ and $\curl W=\Delta u=-\mu u$
			in the weak sense, and $|W|=|\nabla u|$,
			applying \eqref{field-local-identities-div-curl} to $W$ gives \eqref{rotated-derivative} and \eqref{rotated-mass}.
		\end{proof}

	The following lemma shows that the weak tangential boundary
	condition is preserved under multiplication by a function
	that is real on the boundary and addition of a field whose
	coefficient belongs to $H^1_0(\D)$.
	\begin{lemma}\label{admissible}
		Let $W\in\mathscr X$ and let $\beta\in C^\infty(\overline\D;\mathbb C)$
		be real on $\partial\D$. Suppose that $F_Z-\beta F_W\in H^1_0(\D)$
		and that $Z,\diver Z,\curl Z$ are square integrable. Then $Z\in\mathscr X$.
	\end{lemma}

	\begin{proof}
    Let $Y_0$ be the tangent field on $\Omega$ such that $F_{Y_0}=F_Z-\beta F_W$.
		Since the conformal map induces an isomorphism between $H^1_0(\D)$ and $H^1_0(\Omega)$ by pullback,
		we have
		$Y_0\in H^1_0(\Omega;T\Omega)$.
		Approximating $Y_0$ by compactly supported smooth fields and
		passing to the limit in \eqref{tang}, we obtain
		$\mathcal T(Y_0,\varphi)=0$ for every $\varphi\in H^1(\Omega)$.
		
		Write $\beta\circ\Psi^{-1}=a+ib$. The conformal correspondence
		gives $a,b\in H^1(\Omega)\cap L^\infty(\Omega)$, and
		$b\in H^1_0(\Omega)$ because $\beta$ is real on $\partial\D$.
		Since multiplication by $i$ represents $J$, the field with
		coefficient $\beta F_W$ is
		\[
		G:=aW+bJW=Z-Y_0.
		\]
		By the assumptions on $Z$ and  $Y_0\in H^1_0(\Omega)$,
		we have $G,\diver G,\curl G\in L^2(\Omega)$.
		It remains to show that $G$ satisfies the weak tangential
		boundary condition.
		
		For $\psi\in C_0^\infty(\Omega)$, we have
		$a\psi\in H^1_0(\Omega)$, by the definition of the weak
		curl of $W$, together with the product rule, gives
		\[
		\begin{aligned}
			\langle\curl(aW),\psi\rangle
			&=-\int_\Omega a\ip{W}{J\nabla\psi} dA
			=-\int_\Omega\ip{W}{J\nabla(a\psi)} dA
			+\int_\Omega\psi\ip{W}{J\nabla a} dA\\
			&=\int_\Omega\bigl(a\curl W+\ip{W}{J\nabla a}\bigr)\psi dA.
		\end{aligned}
		\]
		Hence
		\begin{equation}\label{curl-product}
			\curl(aW)=a\curl W+\ip{W}{J\nabla a}
			\qquad\text{in weak sense}.
		\end{equation}
		Similarly, $\curl(bJW)=b\curl(JW)+\ip{JW}{J\nabla b}$
		in weak sense. Therefore, we have
		\begin{equation}\label{curl-product-test}
			\begin{aligned}
				\mathcal T(G,\varphi)
				&=\int_\Omega\bigl((\curl W)a\varphi
				+\ip{W}{J\nabla(a\varphi)}\bigr) dA\\
				&\quad+\int_\Omega\bigl((\curl(JW))b\varphi
				+\ip{JW}{J\nabla(b\varphi)}\bigr) dA\\
				&=\mathcal T(W,a\varphi)+\mathcal T(JW,b\varphi).
			\end{aligned}
		\end{equation}
		Here $a\varphi\in H^1(\Omega)$, so the term involving $W$
		vanishes because $W\in\mathscr X$. Also,
		$b\varphi\in H^1_0(\Omega)$ and
		$\curl(JW)=\diver W\in L^2(\Omega)$, so the term involving $JW$
		vanishes by the definition of its weak curl. Thus
		\begin{equation}\label{productgreen}
			\mathcal T(G,\varphi)=0.
		\end{equation}
		Adding the identity for $Y_0$ gives $\mathcal T(Z,\varphi)=0$,
		and hence $Z\in\mathscr X$.
	\end{proof}

	\section{The strict two-multiplier inequality on the disk}\label{disksec}
	Let $\mathcal H$ be a separable complex Hilbert space of finite or countably infinite dimension, and put $\mathcal L=L^2(\D;\mathcal H)$. Fix an orthonormal basis $\{e_j\}$ of $\mathcal H$. In what follows, differential operators on $\mathcal H$-valued functions are understood to act componentwise. For $X(z)=\sum_j X_j(z) e_j$, we write
	\[
	\partial_x X=\sum_j(\partial_x X_j(z))e_j,\quad \partial_y X=\sum_j(\partial_y X_j(z))e_j,
	\text{ and } \Delta X=\sum_j(\Delta X_j(z))e_j,
	\]
	with derivatives understood in the weak sense  whenever necessary.
	If $X\in H^1(\D;\mathcal H)$, we also write
	$\|\nabla X\|^2=\|\partial_x X\|^2+\|\partial_y X\|^2$.
	Clearly, these quantities are independent of the choice of basis.

	For $z\in\mathbb D$, let $V(z)$ be a nonnegative bounded self-adjoint linear operator on $\Hh$.
	Assume that $z\mapsto V(z)\alpha$ is measurable for each fixed $\alpha\in\Hh$,
	and that for almost every $z$, the trace is finite,
	\begin{align}\label{eq:trace-V}
		v(z):=\operatorname{tr} V(z)=\sum_j\langle e_j,V(z)e_j\rangle <\infty,
	\end{align}
	We further assume that $v\in L^r(\mathbb D)$ for some $r>1$ and that there exists $\eta<1$ such that
	\begin{equation}\label{small}
		\int_\mathbb D \langle\xi,V\xi\rangle  \le \eta\int_\mathbb D \|\nabla\xi\|^2 
		\quad \text{ for all }\xi\in H^1_0(\mathbb D;\mathcal H).
	\end{equation}
	Let $S_V(z)=V(z)^{1/2}$ be the unique nonnegative square root of $V(z)$,
	where the powers $V(z)^{1/2}$ are understood in the usual spectral sense that $\langle S_V(z)\alpha,S_V(z)\beta\rangle=\langle \alpha,V(z)\beta\rangle$
	for all $\alpha,\beta\in \mathcal H$. The Cauchy--Schwarz inequality for the form
	$\langle \cdot ,V(z)\cdot\rangle$ gives the pointwise estimate
	\begin{equation}\label{potential-pointwise}
		\|S_V(z)\alpha\|^2=\langle \alpha,V(z)\alpha\rangle
		\le v(z)\|\alpha\|^2 .
	\end{equation}
	For $\Hh$-valued functions with $\pz X,S_VX\in\mathcal L$, we define
	\begin{equation}\label{diskform}
		\mathcal A_V(X)
			=4\|\pz X\|_{\mathcal L}^2-\|S_VX\|_{\mathcal L}^2.
	\end{equation}
	
	Let $f\in H^{1/2}(\partial\D;\Hh)$, and let $U_{\rm a}$ be a fixed
	$\Hh$-valued antiholomorphic function in $\D$ such that
	$U_{\rm a}(0)=0$ and $S_VU_{\rm a}\in\mathcal L$. We consider the
	Dirichlet problem
	\begin{equation}\label{base-problem}
		\begin{cases}
			&(-\Delta-V)U=0 \quad\text{ in }\D,\\
			&(U-U_{\rm a})=f \quad \text{on } \partial\D,\\
			& U-U_{\rm a}\in H^1(\D;\Hh),
		\end{cases}
	\end{equation}
	Set $\gamma_1=\cos\theta$,
	$\gamma_2=\sin\theta$ on $\partial\D$, and define their extensions away from the origin by
	\begin{equation}\label{circle-multipliers}
		\Gamma_1(z)=\frac{z+z^{-1}}2,\qquad
		\Gamma_2(z)=\frac{z-z^{-1}}{2i},\qquad z\ne0.
	\end{equation}
	Clearly, $\Gamma_i=\overline{\Gamma_i}=\gamma_i$ on $\partial\D$.
	Since $U_{\rm a}(0)=0$, write $U_{\rm a}=\bar zU_*$ with $U_*$
	antiholomorphic. Then
	\begin{equation}\label{removable-products}
		\overline{\Gamma_1}U_{\rm a}=\tfrac12(1+\bar z^2)U_*,\qquad
		\overline{\Gamma_2}U_{\rm a}=\tfrac1{2i}(1-\bar z^2)U_*.
	\end{equation}
	Thus both products extend antiholomorphically through the origin.
	We next consider the corresponding comparison problems
	\begin{equation}\label{comparison-problems}
		\begin{cases}
			&(-\Delta-V)U_i=0 \quad\text{ in }\D,\\
			&(U_i-\overline{\Gamma_i}U_{\rm a})=\gamma_if \quad \text{on } \partial\D,\\
			& U_i-\overline{\Gamma_i}U_{\rm a}\in H^1(\D;\Hh),
		\end{cases}
	\end{equation}

	We first prove existence and uniqueness for \eqref{base-problem}. In fact,
	we choose a fixed $F\in H^1(\D;\Hh)$ with $F|_{\partial\D}=f$, and seek a solution of the form
	\[
	U=U_{\rm a}+F+W,\qquad W\in H^1_0(\D;\Hh).
	\]
	For $X\in H^1(\D;\Hh)$ and $\varphi\in H^1_0(\D;\Hh)$, define
	the bounded sesquilinear form
	\begin{equation}\label{sesquilinear}
		a_V(X,\varphi)
		=\int_\D\bigl(\ip{\nabla X}{\nabla\varphi}
		-\ip{VX}{\varphi}\bigr) d z.
	\end{equation}
	By the strict form bound \eqref{small} and
	Poincar\'e's inequality give
	\begin{equation}\label{correction-estimates}
		a_V(X,X)\ge(1-\eta)\|\nabla X\|_{L^2(\D;\Hh)}^2
		\ge c\|X\|_{H^1(\D;\Hh)}^2, \quad\text{ for all } X\in H^1_0(\D;\Hh),
	\end{equation}
	with some $c>0$. Hence the sesquilinear  form is coercive on $H^1_0(\D;\Hh)$.
	Because $U_{\rm a}$ is harmonic, the weak form of
	\eqref{base-problem} is
	\[
	a_V(U-U_{\rm a},\varphi)
	=\int_\D\ip{S_VU_{\rm a}}{S_V\varphi} d z
	\qquad(\varphi\in H^1_0(\D;\Hh)).
	\]
	Substituting $U=U_{\rm a}+F+W$, we obtain the weak equation for the zero-boundary correction $W$,
	\begin{equation}\label{base-correction}
		\begin{aligned}
			a_V(W,\varphi) =\ell_F(\varphi),\quad\text{ where }
			\ell_F(\varphi) =-\int_\D\ip{\nabla F}{\nabla\varphi} d z
			+\int_\D\ip{S_V(U_{\rm a}+F)}{S_V\varphi} d z.
		\end{aligned}
	\end{equation}
	By using the pointwise bound \eqref{potential-pointwise}
	for $V$, together with H\"older's inequality and the Sobolev embedding, we have
	\begin{equation}\label{potential-sobolev}
		\|S_VF\|_{\mathcal L}^2
		\le\int_\D v\|F\|^2 d z
		\le\|v\|_{L^r(\D)}\|F\|_{L^{2r/(r-1)}(\D;\Hh)}^2
		\le C\|v\|_{L^r(\D)}\|F\|_{H^1(\D;\Hh)}^2.
	\end{equation}
	Therefore $S_VF\in\mathcal L$ when $F\in H^1(\D;\Hh)$, and the assumption on $U_{\rm a}$ gives
	$S_V(U_{\rm a}+F)\in\mathcal L$. Thus the Cauchy--Schwarz inequality and \eqref{small} gives
	\[
	|\ell_F(\varphi)|
	\le\bigl(\|\nabla F\|_{L^2(\D;\Hh)}
	+\sqrt\eta\,\|S_V(U_{\rm a}+F)\|_{L^2(\D;\Hh)}\bigr)
	\|\nabla\varphi\|_{L^2(\D;\Hh)}\le C\|\varphi\|_{H^1(\D;\Hh)} .
	\]
	Thus $\ell_F$  is bounded.
	The Lax--Milgram theorem therefore gives a unique
	$W\in H^1_0(\D;\Hh)$ satisfying \eqref{base-correction}.
	Consequently $U=U_{\rm a}+F+W$ is a weak solution of
	\eqref{base-problem}. If two solutions have the same boundary data,
	their difference $R$ belongs to $H^1_0(\D;\Hh)$ and satisfies
	$a_V(R,R)=0$. Coercivity gives $R=0$, proving uniqueness.
	
	The same argument applies to \eqref{comparison-problems} once
	$S_V(\overline{\Gamma_i}U_{\mathrm{a}})\in\mathcal L$ has been verified.
	Fix $0<t<1$. The expressions in \eqref{removable-products} are bounded
	on $|z|\le t$, whereas $\Gamma_i$ is bounded on $t<|z|<1$.
	Using the pointwise bound \eqref{potential-pointwise} for $V$ on the inner disk and the assumption
	$S_VU_{\rm a}\in\mathcal L$ on the outer annulus, we obtain
	\begin{equation}\label{singular-mass}
		\|S_V(\overline{\Gamma_i}U_{\rm a})\|_{\mathcal L}^2
		\le C_{t,U_*}\int_{|z|\le t}v d z
		+C_t\int_{t<|z|<1}\|S_VU_{\rm a}\|^2 d z<\infty.
	\end{equation}
	Also, $\gamma_i f\in H^{1/2}(\partial\D;\Hh)$ because $\gamma_i$
	is smooth. Choosing an $H^1$ extension $F_i$ of $\gamma_i f$ and
	repeating the preceding construction gives a unique weak solution
	$U_i=\overline{\Gamma_i}U_{\rm a}+F_i+W_i$, with
	$W_i\in H^1_0(\D;\Hh)$. 
	Thus, for all three solutions, the solution decomposes into the sum of an antiholomorphic part and an $H^1$ part, where the $\pz$ derivatives of the antiholomorphic part vanish.
	Consequently,
	\begin{equation}\label{solution-energy-data}
		\begin{aligned}
			\pz U&=\pz(F+W)\in\mathcal L,
			&\pz U_i&=\pz(F_i+W_i)\in\mathcal L,\\
			S_VU&=S_VU_{\rm a}+S_V(F+W)\in\mathcal L,
			&S_VU_i&=S_V(\overline{\Gamma_i}U_{\rm a})
			+S_V(F_i+W_i)\in\mathcal L.
		\end{aligned}
	\end{equation}
	Hence their $\mathcal A_V$ energies are finite.
	
	The condition $U(0)=0$ supplies the additional constraint in the
	comparison of $U$ with $U_1$ and $U_2$. The following theorem gives
	the resulting inequality for their disk energies.
	
	\begin{theorem}\label{strict}
		Let $U$ and $U_i$, $i=1,2$, be the unique weak solutions of
		\eqref{base-problem} and \eqref{comparison-problems}, respectively.
		If $U(0)=0$, then
		\begin{equation}\label{negative}
			Q_V:=\mathcal A_V(U_1)+\mathcal A_V(U_2)-\mathcal A_V(U)\le0.
		\end{equation}
		The inequality is strict if $\|S_VU\|_{\mathcal L}>0$.
	\end{theorem}

	For the subsequent proof, we replace the arbitrary $H^1$ terms $F$ and $F_i$ by the harmonic extensions $h_{\mathrm{r}}$ and $v_{i,\mathrm{r}}$ of $f$ and $\gamma_i f$, respectively. Here the harmonic extension means the solution of \eqref{base-problem} with $V=0$ and $U_{\mathrm{a}}=0$. The existence and uniqueness for \eqref{base-problem} guarantee that the harmonic extension exists, while regularity theory further shows that it is smooth in the interior and belongs to $H^1$ globally for general Lipschitz domains. Together with the antiholomorphic parts, these give
	\begin{equation}\label{harmdata}
		h=U_{\mathrm{a}}+h_{\mathrm{r}},\qquad
		v_i=\overline{\Gamma_i}U_{\mathrm{a}}+v_{i,\mathrm{r}}.
	\end{equation}
	Thus $h$ and $v_i$ are harmonic in $\D$, and each solution 
	decomposes into the sum of an harmonic part and an $H^1_0$ part
	\begin{equation}\label{harmonic-corrections}
		\begin{aligned}
			U&=h+q_h, &q_h&\in H^1_0(\D;\Hh),\\
			U_i&=v_i+q_{v_i},
			&q_{v_i}&\in H^1_0(\D;\Hh).
		\end{aligned}
	\end{equation}
	Since $U=h+q_h$ and $h$ is harmonic, the zero boundary correction term
	$q_h$ satisfies
	\begin{equation}\label{correct}
		a_V(q_h,\varphi)
		=\int_\D\ip{S_Vh}{S_V\varphi} d z
		\quad\text{ for all }\varphi\in H^1_0(\D;\Hh).
	\end{equation}
	The same equation holds for $q_{v_i}$ with $v_i$ in place of $h$.

	The next lemma gives the boundary relation needed when the disk
	solutions are used to construct comparison fields.
	
	\begin{lemma}\label{boundary-multiplier}
		For the solutions $U,U_i$ of \eqref{base-problem} and
		\eqref{comparison-problems}, there are smooth functions
		$\beta_i$ on $\overline\D$, real on $\partial\D$, such that
		\begin{equation}\label{tracecorrection}
			(\beta_i)|_{\partial\D}=\gamma_i,\qquad
			U_i-\beta_iU\in H^1_0(\D;\Hh),\qquad i=1,2.
		\end{equation}
	\end{lemma}
	\begin{proof}
		Choose a real $\chi\in C^\infty(\overline\D)$ equal to zero
		near $0$ and to one near $\partial\D$, and set
		\begin{equation}\label{cutoff-mult}
			\beta_i=\chi\overline{\Gamma_i},\qquad \beta_i(0)=0.
		\end{equation}
		The cutoff removes the pole of $\Gamma_i$, while
		$\Gamma_i=\gamma_i\in\mathbb R$ on $\partial\D$.
		Using \eqref{harmonic-corrections} and \eqref{harmdata}, we have
		\[
		U_i-\beta_iU
		=(1-\chi)\overline{\Gamma_i}U_{\rm a}
		+(v_{i,\rm r}-\beta_i h_{\rm r})
		+(q_{v_i}-\beta_i q_h).
		\]
		The first term is smooth and supported away from $\partial\D$
		by \eqref{removable-products}. The second belongs to $H^1$ and
		has zero boundary value, since both $v_{i,\rm r}$ and
		$\beta_i h_{\rm r}$ have boundary value $\gamma_i f$.
		The last belongs to $H^1_0$ because $\beta_i$ is smooth.
		Thus all three terms belong to $H^1_0(\D;\Hh)$, proving
		\eqref{tracecorrection}.
	\end{proof}

	\subsection{The harmonic extensions}
	
		Because the antiholomorphic parts of $h$ and $v_i$ have zero $\pz$
		derivative, their derivative energies come from the $H^1$ harmonic
		extensions of the boundary data. The next lemma computes the
		difference of these energies and shows that it depends only on $h(0)$.

	\begin{lemma}\label{harmonic}
		Let $c=h(0)$, and set 
		\begin{equation}\label{disk-weights}
			d_1=\frac{2x}{1+|z|^2},\qquad d_2=\frac{2y}{1+|z|^2},\qquad
			s=\frac{1-|z|^2}{1+|z|^2},\qquad c_*=\frac1{2\pi}.
		\end{equation}
		Then $|d|^2+s^2=1$, and
		\begin{subequations}\label{harmpair}
			\begin{align}
				&h=d_1v_1+d_2v_2+sc,\label{harm-linear}\\
				&4\sum_{i=1}^2\|\pz v_i\|_{\mathcal L}^2-4\|\pz h\|_{\mathcal L}^2
				=c_*^{-1}\|c\|^2 .\label{harm-energy}
			\end{align}
		\end{subequations}
	\end{lemma}
	\begin{proof}
		We first consider a finite Fourier sum
		$f(\theta)=\sum_n f_ne^{in\theta}$, where $f_n\in \Hh$. Since the harmonic extension
		of $e^{in\theta}$ is $r^{|n|}e^{in\theta}$,  we have
		\begin{equation}\label{fourier-shifts1}
			h_{\rm r}(z)=\sum_{n} r^{|n|}f_n e^{in\theta}.
		\end{equation}
		Since $\gamma_1\pm i\gamma_2=e^{\pm i\theta}$, the function
		$v_{1,\rm r}\pm iv_{2,\rm r}$ is the harmonic extension
		of $e^{\pm i\theta}f$, thus
		\begin{align}\label{fourier-shifts2}
			(v_{1,\rm r}\pm iv_{2,\rm r})(z)
			=\sum_n f_n r^{|n\pm1|}e^{i(n\pm1)\theta}. 
		\end{align}
		Multiplying by $\bar z$ and
		$z$, respectively, and using the definitions of $d_1,d_2$, we obtain
		\begin{equation*}
			\begin{aligned}
				d_1v_{1,\rm r}+d_2v_{2,\rm r}
				&=\frac{\bar z(v_{1,\rm r}+iv_{2,\rm r})+z(v_{1,\rm r}-iv_{2,\rm r})}{1+r^2}\\
				&=\sum_n f_ne^{in\theta}
				\frac{r^{1+|n+1|}+r^{1+|n-1|}}{1+r^2}.
			\end{aligned}
		\end{equation*}
		For $n\ne0$ the coefficient in the last sum is $r^{|n|}$, while
		for $n=0$ it is $2r^2/(1+r^2)=1-s$, so the sum equals
		$h_{\rm r}-sf_0$. For the antiholomorphic part, direct calculation gives
		\begin{equation*}
			\sum_i d_i\overline{\Gamma_i}
			=\tfrac12\bigl((d_1+id_2)\bar z+(d_1-id_2)\bar z^{-1}\bigr)=1 .
		\end{equation*}
		Multiplying this identity by $U_{\rm a}$ gives
		$d_1(\overline{\Gamma_1}U_{\rm a})+d_2(\overline{\Gamma_2}U_{\rm a})=U_{\rm a}$,
		since $d_i(0)=U_{\rm a}(0)=0$, the
		identity holds there as well. Adding the regular part computed above,
		we obtain
		\[
		d_1v_1+d_2v_2=U_{\rm a}+h_{\rm r}-sf_0=h-sf_0.
		\]
		This proves \eqref{harm-linear}.
		The identity $|d|^2+s^2=1$ follows directly from the definitions.

			For $n\ge1$, the $n$-th terms in \eqref{fourier-shifts1} is
			$r^ne^{in\theta}f_n=z^nf_n$, and direct integration gives
			\begin{align}
				4\int_\D\|\pz(z^nf_n)\|^2 d z
				=8\pi n^2\|f_n\|^2\int_0^1r^{2n-1} d r
				=4\pi n\|f_n\|^2.
			\end{align}
			The term for $n=0$ is constant, while those for $n<0$ are
			$\bar z^{-n}f_n$, so their $\pz$ derivatives vanish. Since the
			derivatives of terms are orthogonal in $L^2(\D;\Hh)$
			and $\pz U_{\rm a}=0$, we obtain
			\begin{align}
				4\int_\D \|\pz h\|^2dz
				=4\int_\D \|\pz h_{\rm r}\|^2dz
				=4\pi\sum_{n}n_+\|f_n\|^2.
			\end{align}
			where  $n_+=\max(n,0)$,
			Applying the same calculation to series
			\eqref{fourier-shifts2} yields
			\begin{align}
				4\int_\D \|\pz(v_{1,\rm r}\pm i v_{2,\rm r})\|^2dz
				=4\pi\sum_{n }(n\pm1)_+\|f_n\|^2.
			\end{align}
			Hence the parallelogram identity gives
			\begin{align}
				4\int_\D \big( \|\pz v_1\|^2 + \|\pz v_2\|^2 \big) dz
				&=2 \int_\D 
				\|\pz(v_{1,\rm r}+i v_{2,\rm r})\|^2 dz+  2 \int_\D \|\pz(v_{1,\rm r}- i v_{2,\rm r})\|^2 dz\\
				&=2\pi\sum_{n\in\mathbb Z}
				\bigl((n+1)_++(n-1)_+\bigr)\|f_n\|^2.
			\end{align}
			Consequently, we have 
			\begin{align}
				4\sum_{i=1}^2\|\pz v_i\|_{\mathcal L}^2
				-4\|\pz h\|_{\mathcal L}^2
				=2\pi\sum_{n\in\mathbb Z}
				\bigl((n+1)_++(n-1)_+-2n_+\bigr)\|f_n\|^2
				=2\pi\|f_0\|^2.
			\end{align}
			This proves \eqref{harm-energy} for finite Fourier sums.
			For general $f\in H^{1/2}$, Fourier partial sums converge in
			$H^{1/2}$, and their harmonic extensions converge in $H^1$ and
			locally uniformly, so both identities follow by passage to the limit.
	\end{proof}
	
	\subsection{The Green operator and the value at the center}
	Let $A=-\Delta_D$ be the negative Laplacian on $\D$ with 
	Dirichlet boundary conditions, acting componentwise
	on $\mathcal L$, and set $G=A^{-1}$ is the Dirichlet Green operator. 
	The priori estimates implies that  $G:\mathcal L\to H^1_0(\D;\Hh)$ is bounded,
	and is positive, self-adjoint operator on $\mathcal L$,
	thus $A^{-1/2}:\mathcal L\to H^1_0(\D;\Hh)$ is well defined in the usual spectral sense,
	$A^{-1/2}$ is bounded and positive, self-adjoint operator on $\mathcal L$.
	
	Define $Q=S_VA^{-1/2}$.
	For $g\in\mathcal L$, set $\xi=A^{-1/2}g$. Then
	$\xi\in H^1_0(\D;\mathcal H)$ and
	$\|\nabla\xi\|_{L^2(\D;\mathcal H)}=\|g\|_{\mathcal L}$.
	The assumption
	\eqref{small}, with constant $\eta<1$, gives
	\begin{align}\label{eq:Qg-normal}
		\|Qg\|_{\mathcal L}^2
		=\int_\D\ip{\xi}{V\xi} dz
		\le\eta\|\nabla\xi\|_{L^2(\D;\mathcal H)}^2
		=\eta\|g\|_{\mathcal L}^2.
	\end{align}
	Set $T=QQ^*$, then $T=S_VGS_V$. Since \eqref{eq:Qg-normal} gives $\|Q\|\le\sqrt{\eta}$, we have
	\[
	0\le\ip{g}{Tg}=\|Q^*g\|_{\mathcal L}^2
	\le\eta\|g\|_{\mathcal L}^2
	\quad \text{ for all }g\in\mathcal L.
	\]
	Thus $0\le T\le\eta I$. In particular, let $B=I-T$, then $I\ge B\ge(1-\eta)I$ is
	invertible, and we set $R=B^{-1}$. Clearly, $B$ and $R$ are bounded self-adjoint operators on $\mathcal L$.
	
	%
	%

	On the disk, the Dirichlet Green kernel and its
	value with one argument at the center are
	\begin{equation}\label{green-function}
		G(z,t)=c_*\log\left|\frac{1-z\bar t}{z-t}\right|,\qquad
		\ell(z)=G(z,0)=c_*\log(1/|z|), \quad \text{ where } c_*=\frac{1}{2\pi}.
	\end{equation}
	Define $\mathsf E:\mathcal H\to\mathcal L$
	by $(\mathsf E\alpha)(z)=\ell(z)S_V(z)\alpha$.
	By the pointwise estimate $\ip{\alpha}{V(z)\alpha}\le v(z)\|\alpha\|^2$ and $S_V(z)^2=V(z)$ ,
	then H\"older's inequality gives
	\begin{equation}\label{evaluation-bound}
		\|\mathsf E\alpha\|_{\mathcal L}^2
		=\int_\D \ell(z)^2\ip{\alpha}{V(z)\alpha} dz
		\le \|\alpha\|^2\int_\D \ell^2v
		\le \|\alpha\|^2\|\ell\|_{L^{2r'}(\D)}^2\|v\|_{L^r(\D)},
	\end{equation}
	where $r'=r/(r-1)$.
	Since $\ell(z)=c_*\log(1/|z|)$, polar integration gives
	\[
	\int_\D|\ell(z)|^p dz
	=2\pi c_*^p\int_0^1(-\log t)^p t dt
	=2\pi c_*^p\int_0^\infty s^p e^{-2s} ds<\infty
	\quad \text{ for }1\le p<\infty.
	\]
	In particular, $\ell\in L^{2r'}(\D)$. Therefore \eqref{evaluation-bound} and $v\in L^r(\D)$ impliy the boundedness of $E$.

	\begin{lemma}\label{green}
		Let $(X,\psi)$ be any of $(U,h)$, $(U_1,v_1)$, and $(U_2,v_2)$,
		and put $w=S_VX$ and $\widetilde\psi=S_V\psi$. Then $X$ is
		continuous in $\D$, and
		\begin{subequations}\label{three}
			\begin{align}
				w&=R\widetilde\psi,\label{green-resolvent}\\
				\mathcal A_V(X)&=4\|\pz\psi\|_{\mathcal L}^2
				-\ip{\widetilde\psi}{R\widetilde\psi}_{\mathcal L},
				\label{green-value-energy}\\
				X(0)&=\psi(0)+\mathsf E^*w .\label{green-center}
			\end{align}
		\end{subequations}
	\end{lemma}
	\begin{proof}
		We prove the assertions for $(X,\psi)=(U,h)$. Write
		$q_h=U-h\in H^1_0(\D;\Hh)$ and $w=S_VU$. Since
		$w=S_V(h+q_h)$, equation \eqref{correct} takes the form
		\begin{equation}\label{green-correction-weak}
			\int_\D\ip{\nabla q_h}{\nabla\varphi} d z
			=\int_\D\ip{w}{S_V\varphi} d z
			\quad \text{ for all }\varphi\in H^1_0(\D;\Hh).
		\end{equation}
		Taking $\varphi=A^{-1/2}g\in H^1_0(\D;\Hh)$ for any $g\in\mathcal L$ and using
		$Q=S_VA^{-1/2}$, we obtain
		\[
		\int_\D\ip{A^{1/2}q_h}{g} d z
		=\int_\D\ip{w}{Qg} d z
		=\int_\D\ip{Q^*w}{g} d z.
		\]
		Thus $A^{1/2}q_h=Q^*w$. Since $Q=S_VA^{-1/2}$ and
		$w=S_V(h+q_h)$, it follows that
		\begin{equation}\label{green-correction}
			q_h=A^{-1/2}Q^*w,\qquad S_Vq_h=Tw,\qquad
			(I-T)w=S_Vh .
		\end{equation}
		Inverting $I-T$ proves \eqref{green-resolvent} for $U$.

		Taking $\varphi\in C_0^\infty(\D;\mathcal H)$,
		integration by parts, together with the fact that $h$ is harmonic in the interior, gives
		$4\int_\D\ip{\pz\varphi}{\pz h} d z
		=-\int_\D\ip{\varphi}{\Delta h} d z=0$.
		Since $\pz h\in\mathcal L$, the left-hand side is continuous in
		the $H^1_0$ norm, so the identity extends by density to every
		$\varphi\in H^1_0(\D;\Hh)$. Integration by parts also gives
		\[
		4\int_\D\|\pz\varphi\|^2 d z
		=-\int_\D\ip{\varphi}{\Delta\varphi} d z
		=\int_\D\|\nabla\varphi\|^2 d z.
		\]
		By density, this identity holds for $q_h\in H^1_0(\D;\Hh)$.
		Testing \eqref{green-correction-weak} with $q_h$ now yields
		\begin{equation}\label{green-test}
			\int_\D\ip{\pz q_h}{\pz h} d z=0,\qquad
			4\|\pz q_h\|_2^2=\|\nabla q_h\|_2^2
			=\int_\D\ip{S_Vq_h}{w} d z.
		\end{equation}
		Expanding the energy term of $U=h+q_h$, together with \eqref{green-test} and $w=S_V(h+q_h)$, we obtain
		\begin{equation*}
			\begin{aligned}
				\mathcal A_V(U)
				&=4\|\pz h\|_{\mathcal L}^2+\|\nabla q_h\|_2^2-\|S_V(h+q_h)\|_{\mathcal L}^2\\
				&=4\|\pz h\|_{\mathcal L}^2
				+\int_\D\ip{S_Vq_h}{w} d z-\int_\D\langle S_V(h+q_h), w\rangle dz \\
				&=4\|\pz h\|_{\mathcal L}^2-\int_\D\ip{S_Vh}{w} d z.
			\end{aligned}
		\end{equation*}
		Together with \eqref{green-resolvent}, this proves the energy
		identity \eqref{green-value-energy}.

		To evaluate the solution at
		the center, set $p=2r/(r+1)>1$ and let $p'$ be its conjugate
		exponent. 
		By using the pointwise bound \eqref{potential-pointwise}
		for $V$, together with H\"older's inequality, we have
		\begin{equation}\label{green-source}
			\|S_Vg\|_{L^p(\D;\Hh)}
			\le\|  v\|_{L^{ r}(\D)}^{1/2}\|g\|_{\mathcal L}
			\qquad \text{ for any } g\in\mathcal L.
		\end{equation}
		In particular, $S_Vw\in L^p(\D;\Hh)$ since $w\in \mathcal L$ by \eqref{solution-energy-data}.
		Since $0\le G(z,t)\le c_*\log(2/|z-t|)$, its kernel satisfies bound
		$\sup_{z\in\D}\|G(z,\cdot)\|_{p'}<\infty$.
		Approximating $S_Vw$ in $L^p$ by smooth compactly function, their
		Dirichlet solutions converge in $H^1_0$, and their Green integrals
		converge uniformly by the kernel bound. Hence $q_h$ has
		the Green representation
		\begin{equation}\label{green-integral}
			q_h(z)=\int_\D G(z,t)S_V(t)w(t) d t .
		\end{equation}
		Therefore, for $z_0\in \D$, by using the H\"older's inequality
		and the $L^{p'}$ continuity of $G(z,\cdot)$, we have
		\begin{equation*}
			\|q_h(z)-q_h(z_0)\|
			\le\|G(z,\cdot)-G(z_0,\cdot)\|_{L^{p'}(\D)}
			\|S_Vw\|_{L^p(\D;\Hh)}\xrightarrow{z\to z_0}0 .
		\end{equation*}
		This proves continuity of $q_h$ and hence of $U=h+q_h$ in
		the disk.
		Finally, taking $z=0$ in \eqref{green-integral}, where
		$G(0,t)=\ell(t)$, we have $q_h(0)=\mathsf E^*w$, and \eqref{green-center} follows from $U=h+q_h$.
		The same argument with $(U_i,v_i)$ in place of $(U,h)$ proves
		all three assertions for $i=1,2$.
	\end{proof}

	\subsection{Strict positivity of the kernel}
	Let $G$ be the Dirichlet Green kernel of the disk and define $H$ and
	$J$ for $z\ne t$ by
	\begin{equation}\label{kernel-data}
		H(z,t)=\frac{|z-t|^2}{(1-|z|^2)(1-|t|^2)},\qquad
		k(z,t)=\frac1{1+H(z,t)},\qquad
		J(z,t)=c_*-2H(z,t)G(z,t).
	\end{equation}
	Set $J(z,z)=c_*$ and $k(z,z)=1$ when $z=t$.

	\begin{lemma}\label{poskernel}
		For every nonzero $\Hh$-valued Borel measure $m$ of finite total variation
		on $\D$, we have
		\begin{equation}\label{measure-strict}
			\iint J(z,t)\ip{ d m(z)}{ d m(t)}>0.
		\end{equation}
	\end{lemma}
	\begin{proof}
		The identity
		$|1-z\bar t|^2-|z-t|^2=(1-|z|^2)(1-|t|^2)$ gives
		$G=(c_*/2)\log(1+H^{-1})$. 
		Integrating
		$\tau/(H+\tau)=1-H/(H+\tau)$ over $(0,1)$ and then expanding
		in powers of $k=(1+H)^{-1}$, we obtain
		\begin{equation}\label{Lseries}
			\begin{aligned}
				J/c_*&=1-H\log(1+H^{-1})
				=\int_0^1\frac{\tau}{H+\tau} d\tau\\
				&=\int_0^1\frac{\tau k}{1-(1-\tau)k} d\tau
				=\sum_{m\ge1}\frac{k^m}{m(m+1)} .
			\end{aligned}
		\end{equation}
		The series converges uniformly for $0\le k\le1$.
		Expanding $(1-z\bar t)^{-1}$ as geometric series gives
		\begin{equation}\label{kernel-features}
			k(z,t)
			=\frac{(1-|z|^2)(1-|t|^2)}{|1-z\bar t|^2}
			=\sum_{j,l\ge0}\overline{\phi_{jl}(z)}\,\phi_{jl}(t),\quad
			\text{ where }
			\phi_{jl}(z)=(1-|z|^2)z^j\bar z^l.
		\end{equation}
		Taking $t=z$ in \eqref{kernel-features} gives
		$\sum_{j,l}|\phi_{jl}(z)|^2=1$.
		By Cauchy--Schwarz, the sum of the
		absolute values of the terms in \eqref{kernel-features} is at
		most $1$, together with that $m$ has finite total variation, we may integrate
		term by term and obtain
		\begin{equation}\label{eq:kernel-square-sum}
			\iint k(z,t)\langle dm(z),dm(t)\rangle
			=\sum_{j,l\ge0}
			\left\|\int_\D (1-|z|^2)z^j\bar z^l\,dm(z)\right\|_{\mathcal H}^2.
		\end{equation}

		For each $n\ge1$, multiplying the series \eqref{kernel-features} gives
		\begin{equation}\label{eq:kernel-power-expansion}
			k(z,t)^n
			=\sum_{j_1,l_1,\ldots,j_n,l_n\ge0}
			\overline{\prod_{a=1}^n\phi_{j_a l_a}(z)}
			\prod_{a=1}^n\phi_{j_a l_a}(t).
		\end{equation}
		Since the sum of the absolute values of the terms in \eqref{eq:kernel-power-expansion}
		is at most $1$ and $m$ has finite total variation, we may integrate
		term by term and obtain
		\begin{equation}\label{eq:kernel-power-integral}
			\iint k(z,t)^n\langle dm(z),dm(t)\rangle
			=\sum_{j_1,l_1,\ldots,j_n,l_n\ge0}
			\left\|\int_\D\prod_{a=1}^n\phi_{j_a l_a}(z)\,dm(z)\right\|_{\mathcal H}^2
			\ge0.
		\end{equation}
		Using \eqref{eq:kernel-power-integral} and the uniformly convergent
		series \eqref{Lseries} for $J$, we obtain
		\begin{equation}\label{eq:kernel-J-lower-bound}
			\begin{aligned}
				\iint J(z,t)\langle dm(z),dm(t)\rangle
				&=c_*\sum_{n=1}^\infty\frac{1}{n(n+1)}
				\iint k(z,t)^n\langle dm(z),dm(t)\rangle\\
				&\ge \frac{c_*}{2}\iint k(z,t)\langle dm(z),dm(t)\rangle.
			\end{aligned}
		\end{equation}
		
		Suppose $\iint k(z,t)\langle dm(z),dm(t)\rangle=0$ and set
		$d\widetilde m(z)=(1-|z|^2)\,dm(z)$. Then \eqref{eq:kernel-square-sum} means
		\begin{equation}\label{eq:kernel-weighted-monomials}
			\int_\D z^j\bar z^l\,d\widetilde m(z)=0
			\quad \text{ for all }j,l\ge0.
		\end{equation}
		For $e\in\mathcal H$, define a measure $\mu_e$ on
		$\overline\D$ by $\mu_e(E)=\langle e,\widetilde m(E\cap\D)\rangle$
		for Borel sets $E$. Then we have
		\begin{equation}\label{eq:kernel-polynomial-integral}
			\int_{\overline\D}p(z,\bar z)\,d\mu_e(z)=0
			\qquad\text{for every polynomial }p.
		\end{equation}
		Since polynomials in $z,\bar{z}$ are dense in continuous function, we have
		$\int_{\overline\D}g(z,\bar z) d\mu_e=0$ for any continuous function $g$.
		Thus $\mu_e=0$ for every $e$, which implies $\widetilde m=0$, and hence $m=0$, a contradiction.
		Therefore the right hand side in \eqref{eq:kernel-J-lower-bound} is strictly positive, which proves the claim.
	\end{proof}
	
	\subsection{Operator factorization}
	
	Write $j(z)=J(z,0)$ and define $\mathsf U,\mathsf J_0:\Hh\to\mathcal L$
	and $\mathsf J:\mathcal L\to\mathcal L$ by
	\begin{equation}\label{kernel-operators}
		\begin{aligned}
			(\mathsf U\alpha)(z)&=S_V(z)\alpha,&
			(\mathsf J_0\alpha)(z)&=j(z)S_V(z)\alpha,\\
			(\mathsf Ju)(z)&=S_V(z)\int_\D J(z,t)S_V(t)u(t) d t .
		\end{aligned}
	\end{equation}
	For an operator
	$K:\mathcal H\to\mathcal L$, the Hilbert--Schmidt norm is defined by
	$\|K\|_{\rm HS}^2=\sum_j\|Ke_j\|_{\mathcal L}^2$ for any
	orthonormal basis $(e_j)$ of $\mathcal H$.
	Using $S_V^2=V$ and
	interchanging the nonnegative sums and integrals gives
	\begin{equation}\label{hs-norms}
		\begin{aligned}
			\|\mathsf U\|_{\rm HS}^2
			&=\sum_j\int_\D\|S_V(z)e_j\|_{\mathcal H}^2 dz
			=\int_\D\sum_j\ip{e_j}{V(z)e_j} dz
			=\int_\D v(z) dz,\\
			\|\mathsf E\|_{\rm HS}^2
			&=\sum_j\int_\D\ell(z)^2\|S_V(z)e_j\|_{\mathcal H}^2 dz
			=\int_\D\ell(z)^2v(z) dz,\\
			\|\mathsf J_0\|_{\rm HS}^2
			&=\sum_i\int_\D |j(z)|^2\|S_V(z)e_i\|_{\mathcal H}^2 dz
			=\int_\D |j(z)|^2v(z) dz.
		\end{aligned}
	\end{equation}
	All three integrals are finite since $v,\ell^2v\in L^1(\D)$ and
	$j$ is bounded. Thus, $\mathsf U$, $\mathsf E$ and $\mathsf J_0$ are Hilbert–Schmidt and hence compact.
	Cauchy--Schwarz inequality and pointwise bound \eqref{potential-pointwise}
	therefore gives
	\begin{equation}\label{eq:S-V-estimate}
		\|S_Vu\|_{L^1(\D;\mathcal H)}
		\le \Bigl(\int_\D v(z)dz\Bigr)^{1/2}\|u\|_{\mathcal L} = \|\sqrt v\|_{L^2(\D)}\|u\|_{\mathcal L}.
	\end{equation}
	Using $|J|\le c_*$ and \eqref{eq:S-V-estimate}, we obtain
	\begin{align}\label{kernel-bounds}
		|\ip{u}{\mathsf J y}|
		&\le c_*\|S_Vu\|_{L^1(\D;\Hh)}\|S_Vy\|_{L^1(\D;\Hh)}
		\le c_*\|v\|_{L^1(\D)}\|u\|_{\mathcal L}\|y\|_{\mathcal L}.
	\end{align}
	Thus $\mathsf J$ is bounded.
	Moreover,  the symmetry $J(z,t)=J(t,z)$ gives $\mathsf J^*=\mathsf J$ by interchanging $z$ and $t$.

	Let $\mathsf D_i$ and $S$ be multiplication on $\mathcal L$ by the
	real functions $d_i$ and $s$, respectively, where $d_i$ and $s$ are defined in \eqref{disk-weights}.
	On $\Hh$, we put
	\begin{equation}\label{center-operators}
		\mathsf L=\mathsf E^*R\mathsf E,\qquad
		\mathsf Q=(\mathsf L+c_*I_{\Hh})^{-1}.
	\end{equation}
	The four operators $M,C,P,N$ act on $\mathcal L$,  and $\mathsf K$.
	They are defined by
	\begin{equation}\label{factor-operators}
		\begin{aligned}
			M&=\sum_i\mathsf D_iB\mathsf D_i,& C&=B+S\mathsf U\mathsf E^*,\\
			P&=B+c_*^{-1}\mathsf E\mathsf E^*,& N&=CP^{-1}C^*,\\
			\mathsf K&=B+\mathsf J-(\mathsf E-\mathsf J_0)\mathsf Q
			(\mathsf E-\mathsf J_0)^* .
		\end{aligned}
	\end{equation}
	Here the inverses are well defined because
	$\mathsf L+c_*I_{\Hh}\ge c_*I_{\Hh}$ and
	$P\ge B\ge(1-\eta)I$.
	\Needspace{5\baselineskip}
	\begin{lemma}\label{factor}
		The operators $N$ and $M$ satisfy
		\begin{equation}\label{factorization}
			B-M=S(B+\mathsf J)S-c_*S\mathsf U\mathsf U^*S,\qquad
			N-M=S\mathsf KS.
		\end{equation}
	\end{lemma}
	
	\begin{proof}
		By the definitions of $d=(d_1,d_2),s$ in \eqref{disk-weights}, a direct calculation gives
		\begin{equation}\label{disk-chord}
			\begin{aligned}
				1-d(z)\cdot d(t)-s(z)s(t)
				&=1-\frac{4\operatorname{Re}(z\bar t)+(1-|z|^2)(1-|t|^2)}
				{(1+|z|^2)(1+|t|^2)}\\
				&=\frac{2|z-t|^2}{(1+|z|^2)(1+|t|^2)}
				=2s(z)s(t)H(z,t).
			\end{aligned}
		\end{equation}
		where the last equality uses the definition of $H$ in
		\eqref{kernel-data}.
		For $t=z$, \eqref{disk-chord} also gives $I-\sum_i\mathsf D_i^2=S^2$,
		Multiplying by $G$ and using
		$2HG=c_*-J$ gives the kernel identity
		\begin{equation}\label{kernel-factor}
			\begin{aligned}
				(1-d(z)\cdot d(t))S_V(z)G(z,t)S_V(t)
				=s(z)S_V(z)\bigl(G(z,t)+c_*-J(z,t)\bigr)S_V(t)s(t).
			\end{aligned}
		\end{equation}
		
		For $u,y\in\mathcal L$, the estimate \eqref{green-source}
		gives $S_Vu,S_Vy\in L^p$, apply H\"older inequality with $p=2r/(r+1)>1$, we have
		\begin{equation}
			\begin{aligned}
				\iint_{\D^2}\|S_V(z)u(z)\|G(z,t)\|S_V(t)y(t)\| d t d z
				\le \|S_Vu\|_{L^1}
				\sup_{z\in\D}\|G(z,\cdot)\|_{L^{p'}}
				\|S_Vy\|_{L^p}<\infty.
			\end{aligned}
		\end{equation}
		The remaining constant and $J$ kernels in\eqref{kernel-factor} are also bounded, and their integrals can
		controlled by \eqref{kernel-bounds}, so Fubini's theorem
		applies. using $T=S_VGS_V$ and \eqref{kernel-factor},  we have
		\begin{align}
			\int_\D\ip{u}{(T-\sum_i\mathsf D_iT\mathsf D_i)y} d z
			&=\iint_{\D^2}(1-d(z)\cdot d(t))G(z,t)
			\ip{S_V(z)u(z)}{S_V(t)y(t)} d t d z\\
			&=\iint_{\D^2}s(z)s(t)\bigl(G(z,t)+c_*-J(z,t)\bigr)
			\ip{S_V(z)u(z)}{S_V(t)y(t)} d t d z\\
			&=\int_\D\ip{u}{S(T+c_*\mathsf U\mathsf U^*-\mathsf J)Sy} d z.
		\end{align}
		Since $u,y\in\mathcal L$ are arbitrary, we obtain
		\begin{equation}\label{Tfactor}
			T-\sum_i\mathsf D_iT\mathsf D_i=S(T+c_*\mathsf U\mathsf U^*-\mathsf J)S .
		\end{equation}
		Substituting $B=I-T$ into $M=\sum_i\mathsf D_i B\mathsf D_i$ and expanding, together
		with $I-\sum_i\mathsf D_i^2=S^2$ and \eqref{Tfactor}, we have
		\begin{equation}
			\begin{aligned}\label{eq:factor-B-minus-M}
				B-M=& \Bigl(I-\sum_i\mathsf D_i^2\Bigr)
				-\Bigl(T-\sum_i\mathsf D_iT\mathsf D_i\Bigr)=
				S^2-S(T+c_*\mathsf U\mathsf U^*-\mathsf J)S \\
				=&S(B+\mathsf J)S-c_*S\mathsf U\mathsf U^*S .
			\end{aligned}
		\end{equation}
		This gives the first identity in \eqref{factor-operators}.
		
		We next compute $P^{-1}$. 
		Multiplying the definition $P=B+c_*^{-1}\mathsf E\mathsf E^*$ by $R$ and using $BR=I$, we have
		$PR = I+c_*^{-1}\mathsf E\mathsf E^*R$,
		Multiplying again by $\mathsf E\mathsf Q$ and using the definition $\mathsf L=\mathsf E^*R\mathsf E$, we obtain
		$PR\mathsf E\mathsf Q
		= \bigl(\mathsf E+c_*^{-1}\mathsf E\mathsf L\bigr)\mathsf Q
		= c_*^{-1}\mathsf E(c_*I_{\mathcal H}+\mathsf L)\mathsf Q
		= c_*^{-1}\mathsf E$.
		Therefore, we have
		\begin{align}
			P\bigl(R-R\mathsf E\mathsf Q\mathsf E^*R\bigr)
			= PR-(PR\mathsf E\mathsf Q)\mathsf E^*R
			= I+c_*^{-1}\mathsf E\mathsf E^*R-c_*^{-1}\mathsf E\mathsf E^*R=I.
		\end{align}
		Therefore, we have
		\begin{equation}\label{inverse-P}
			P^{-1}=R-R\mathsf E\mathsf Q\mathsf E^*R.
		\end{equation}
		
		Write $F_0=S\mathsf U$,
		using \eqref{inverse-P},  $C=B+F_0\mathsf E^*$ and the definition $\mathsf L=\mathsf E^*R\mathsf E$,
		together with the self-adjointness of $B$ and $R$,
		we obtain
		\begin{equation}\label{eq:factor-N-expansion-mid}
			\begin{aligned}
				&N=CP^{-1}C^*=CRC^*-(CR\mathsf E)\mathsf Q(CR\mathsf E)^*\\
				&=\big(B+F_0\mathsf E^*+\mathsf EF_0^*+F_0\mathsf LF_0^*\big)-(\mathsf E+F_0\mathsf L)\mathsf Q(\mathsf E^*+\mathsf LF_0^*)
			\end{aligned}
		\end{equation}
		where the last equality uses $CR\mathsf E=\mathsf E+F_0\mathsf L$.
		Since $\mathsf L=\mathsf Q^{-1}-c_*I_{\mathcal H}$, we have
		$\mathsf E+F_0\mathsf L=(\mathsf E-c_*F_0)+F_0\mathsf Q^{-1}$.
		Thus
		\begin{equation*}
			\begin{aligned}
				(\mathsf E+F_0\mathsf L)\mathsf Q(\mathsf E+F_0\mathsf L)^*
				&=(\mathsf E-c_*F_0)\mathsf Q(\mathsf E-c_*F_0)^*
				+(\mathsf E-c_*F_0)F_0^*
				+F_0(\mathsf E-c_*F_0)^*
				+F_0\mathsf Q^{-1}F_0^*\\
				&=(\mathsf E-c_*F_0)\mathsf Q(\mathsf E-c_*F_0)^*
				+\mathsf EF_0^*+F_0\mathsf E^*+F_0\mathsf LF_0^*
				-c_*F_0F_0^*.
			\end{aligned}
		\end{equation*}
		Together with \eqref{eq:factor-N-expansion-mid}, we obtain
		\begin{equation}\label{reduced-N}
			\begin{aligned}
				N-B
				&=\big(F_0\mathsf E^*+\mathsf EF_0^*+F_0\mathsf LF_0^*\big)-(\mathsf E+F_0\mathsf L)\mathsf Q(\mathsf E^*+\mathsf LF_0^*)\\
				&=c_*F_0F_0^*
				-(\mathsf E-c_*F_0)\mathsf Q(\mathsf E-c_*F_0)^*,
			\end{aligned}
		\end{equation}
		Taking $t=0$ in \eqref{disk-chord} gives
		$1-s=2sH(\cdot,0)$ and $2HG=c_*-J$ gives
		$2H(\cdot,0)\ell=c_*-j$, we have
		\begin{equation}\label{eq:factor-center-identity}
			\ell-c_*s=s(\ell-j),
			\qquad
			\mathsf E-c_*F_0=S(\mathsf E-\mathsf J_0).
		\end{equation}
		Substitute \eqref{eq:factor-center-identity} into
		\eqref{reduced-N}, and adding \eqref{reduced-N} with \eqref{eq:factor-B-minus-M},
		the terms $c_*F_0F_0^*$
		cancel, giving $N-M=S\mathsf KS$.
		This gives the second identity in \eqref{factorization} .

	\end{proof}
	
	\begin{lemma}\label{invert}
		For every $0\ne u\in\mathcal L$, $\ip{u}{\mathsf K u}>0$.
		Moreover, $C$ is boundedly invertible on $\mathcal L$.
	\end{lemma}
	\begin{proof}
		Fix $u\in\mathcal L$. 
		Let $\mu_0$ be the Dirac
		measure at $0$ and set
		\begin{equation}\label{eq:factor-alpha-measure}
			\alpha=-\mathsf Q(\mathsf E-\mathsf J_0)^*u,\qquad
			dm(z)=S_V(z)u(z) dz-\alpha d\mu_0(z).
		\end{equation}
		The $L^1$ bound for $S_Vu$ in \eqref{eq:S-V-estimate} ensures that $m$ has finite total variation.
		Since $J(z,0)=J(0,z)=j(z)$ and $J(0,0)=c_*$, we obtain
		\begin{equation}\label{eq:factor-J-measure}
			\begin{aligned}
				\iint J(z,t)\langle dm(z),dm(t)\rangle  =&  \iint J(z,t)
				\langle S_V(z)u(z), S_V(t)u(t) \rangle dz dt\\
				& -\int  J(z,0)\ip{S_V(z)u(z)}{\alpha}  d z
				-\int  J(0,t)\ip{\alpha}{S_V(t)u(t)}  dt
				+J(0,0)\|\alpha\|_{\mathcal H}^2\\
				=&\ip{u}{\mathsf Ju}
				-\bigl(\ip{u}{\mathsf J_0\alpha}
				+\overline{\ip{u}{\mathsf J_0\alpha}}\bigr)
				+c_*\|\alpha\|_{\mathcal H}^2\\
				=&\ip{u}{\mathsf Ju}-2\Re\ip{u}{\mathsf J_0\alpha}
				+c_*\|\alpha\|_{\mathcal H}^2.
			\end{aligned}
		\end{equation}
		On the other hand,
		using $B^{-1}=R$ and $\mathsf E^*R\mathsf E=\mathsf L$ to expand
		the norm, then applying \eqref{eq:factor-J-measure}, we obtain
		\begin{equation}\label{eq:xx1}
			\begin{aligned}
				\|B^{1/2}u+B^{-1/2}\mathsf E\alpha\|_{\mathcal L}^2
				&=\ip{u}{Bu}+2\Re\ip{u}{\mathsf E\alpha}
				+\ip{\alpha}{\mathsf L\alpha}\\
				&=\ip{u}{(B+\mathsf J)u}
				+2\Re\ip{u}{(\mathsf E-\mathsf J_0)\alpha}\\
				&\quad+\ip{\alpha}{(\mathsf L+c_*I_{\mathcal H})\alpha}
				-\iint J(z,t)\langle dm(z),dm(t)\rangle.
			\end{aligned}
		\end{equation}
		Since
		$(\mathsf L+c_*I_{\mathcal H})\alpha
		=-(\mathsf L+c_*I_{\mathcal H})\mathsf Q(\mathsf E-\mathsf J_0)^*u
		=-(\mathsf E-\mathsf J_0)^*u$, we have
		\begin{equation}\label{eq:xx2}
			\begin{aligned}
				2\Re\ip{u}{(\mathsf E-\mathsf J_0)\alpha}
				+\ip{\alpha}{(\mathsf L+c_*I_{\mathcal H})\alpha}
				=\ip{u}{(\mathsf E-\mathsf J_0)\alpha}
				=-\ip{u}{(\mathsf E-\mathsf J_0)\mathsf Q(\mathsf E-\mathsf J_0)^*u}.
			\end{aligned}
		\end{equation}
		Substituting \eqref{eq:xx2} into \eqref{eq:xx1} and using the definition of $\mathsf K$ in \eqref{factor-operators}, we have
		\begin{equation}\label{eq:factor-K-squares}
			\ip{u}{\mathsf Ku}
			=\|B^{1/2}u+B^{-1/2}\mathsf E\alpha\|_{\mathcal L}^2
			+\iint J(z,t)\langle dm(z),dm(t)\rangle.
		\end{equation}
		Both terms in \eqref{eq:factor-K-squares} are nonnegative by \eqref{measure-strict}.
		If their sum is zero, then
		$m=0$. By \eqref{eq:factor-alpha-measure}, $m(\{0\})=-\alpha$, so
		$\alpha=0$. The first term in \eqref{eq:factor-K-squares} becomes
		\begin{align}
			\|B^{1/2}u\|_{\mathcal L}^2= \langle u,Bu\rangle_{\mathcal L}\ge(1-\eta)\|u\|_{\mathcal L}^2,
		\end{align}
		where the last inequality follows from $B\ge(1-\eta)I$, thus we have $u=0$.
		This proves that $K$ is positive operator.

		Finally, if $C^*y=0$, the definition \eqref{factor-operators} implies $Ny=0$,
		then \eqref{factorization} and the positivity of $K$ gives
		\begin{equation}\label{eq:factor-C-injective}
			0=\ip{y}{My}+\ip{Sy}{\mathsf K Sy} \ge \ip{y}{My}
			\ge(1-\eta)\int_\D |d(z)|^2\|y(z)\|_{\mathcal H}^2dz.
		\end{equation}
		Here we used $M=\sum_i\mathsf D_iB\mathsf D_i$ and $B\ge(1-\eta)I$.
		Since $|d(z)|>0$ for $z\ne0$, \eqref{eq:factor-C-injective} gives
		$y=0$, thus cokernel is $\ker C^*=0$.
		\eqref{hs-norms} implies $\mathsf U$ and $\mathsf E$ are compact,
		so $C-B=S\mathsf U\mathsf E^*$
		is compact. Thus $C$ is Fredholm of index zero, and $\ker C=\{0\}$.
		Therefore $C$ is bijective and its inverse is bounded.
	\end{proof}
	
	\subsection{Proof of the two-multiplier inequality}
	
	\begin{proof}[Proof of \Cref{strict}]
		Put $w=S_VU$, $\widetilde h=S_Vh$ and
		$\widetilde v_i=S_Vv_i$. Applying \Cref{green} to the pair $(U,h)$, 
		together with $R=B^{-1}$ and the assumption $U(0)=0$, we have
		\begin{equation}\label{center-constraint}
			\widetilde h=Bw,\qquad c=-\mathsf E^*w,
		\end{equation}
		the self-adjointness of $B$ therefore gives $\ip{\widetilde h}{R\widetilde h}=\ip{Bw}{w}= \ip{w}{Bw}$.
		Thus the definition of $P$ gives
		\begin{align}\label{eq:aab1}
			c_*^{-1}\|c\|^2+\ip{\widetilde h}{R\widetilde h}=c_*^{-1}\ip{\mathsf E^*w}{\mathsf E^*w}+\ip{w}{Bw}
			=\ip{w}{(B+c_*^{-1}\mathsf E\mathsf E^*)w}=\ip{w}{Pw}.
		\end{align}
		Applying the identity \eqref{green-value-energy}
		to $h,v_1,v_2$, together with \eqref{harm-energy} and \eqref{eq:aab1}, we obtain
		\begin{equation}\label{form-P}
			\begin{aligned}
				Q_V&= 4\sum_i\|\partial_zv_i\|_{\mathcal L}^2-\sum_i\langle\widetilde v_i,R\widetilde v_i\rangle
				-4\|\partial_zh\|_{\mathcal L}^2
				+\langle\widetilde h,R\widetilde h\rangle
				\\
				&=c_*^{-1}\|c\|^2+\ip{\widetilde h}{R\widetilde h}
				-\sum_i\ip{\widetilde v_i}{R\widetilde v_i}\\
				&
				=\ip{w}{Pw}-\sum_i\ip{\widetilde v_i}{R\widetilde v_i}.
			\end{aligned}
		\end{equation}

		Multiplying the identity \eqref{harm-linear} by $S_V$ gives
		\begin{align}
			S_V(z)h(z)=\sum_i d_i(z)S_V(z)v_i(z)+s(z)S_V(z)c,
		\end{align}
		so $\widetilde h=\sum_i\mathsf D_i\widetilde v_i+S\mathsf Uc$.
		Using the definition of $C$ in \eqref{factor-operators} and \eqref{center-constraint}, we obtain
		\begin{align} \label{harmonic-constraint}
			\sum_i\mathsf D_i\widetilde v_i=Bw-S\mathsf Uc=Cw.
		\end{align}

		Since $C^*$ is invertible by  \Cref{invert}, we may define
		\begin{equation}\label{disk-residuals}
			y=(C^*)^{-1}Pw,\qquad r_i=\widetilde v_i-B\mathsf D_i y .
		\end{equation}
		The relation $C^*y=Pw$ and the definition of $N$, together
		with the self-adjointness of $P$, we have
		\begin{align}\label{eq:xxb2}
			Ny=CP^{-1}C^*y=Cw,\quad
			\ip{y}{Cw}=\ip{C^*y}{w}=\ip{Pw}{w}=\ip{w}{Pw}.
		\end{align}
		Using $N-M=S\mathsf KS$ from \eqref{factorization} and \eqref{eq:xxb2}, we obtain
		\[
		\ip{y}{My}=\ip{y}{Ny}-\ip{Sy}{\mathsf KSy}
		=\ip{w}{Pw}-\ip{Sy}{\mathsf KSy}.
		\]
		Since $R=B^{-1}$ is self-adjoint, expanding the expression for each $r_i$ gives
		\begin{align}
			\ip{r_i}{Rr_i}
			=\ip{\widetilde v_i}{R\widetilde v_i}
			-2\Re\ip{\mathsf D_i y}{\widetilde v_i}
			+\ip{\mathsf D_i y}{B\mathsf D_i y}.
		\end{align}
		Summing over $i$,
		using \eqref{harmonic-constraint} and the definition of $M$ in \eqref{factor-operators},
		and noting that self-adjointness in \eqref{eq:xxb2} implies that $\langle w, Pw\rangle$ is real, we obtain
		\begin{equation}\label{residual-sum}
			\begin{aligned}
				\sum_i\ip{r_i}{Rr_i}
				&=\sum_i\ip{\widetilde v_i}{R\widetilde v_i}
				-2\Re\ip{y}{Cw}+\ip{y}{My}\\
				&=\sum_i\ip{\widetilde v_i}{R\widetilde v_i}
				-\ip{w}{Pw}-\ip{Sy}{\mathsf KSy}.
			\end{aligned}
		\end{equation}

		Eliminating $\sum_i\ip{\widetilde v_i}{R\widetilde v_i}$ between
		\eqref{form-P} and \eqref{residual-sum} yields
		\begin{equation}\label{exactsquares}
			Q_V=-\sum_i\ip{r_i}{Rr_i}-\ip{Sy}{\mathsf K Sy}.
		\end{equation}
		Since $R$ is positive, and $\mathsf K$ is strictly positive
		by \Cref{invert}, so 
		$Q_V\le0$. If $\|S_VU\|_{\mathcal L}>0$, then $w\ne0$,
		$P$ is invertible, so $C^*y=Pw$ implies
		$y\ne0$.
		Therefore $Sy\ne0$. Strict positivity of $\mathsf K$
		makes \eqref{exactsquares} negative, proving the strict assertion.
	\end{proof}

	\section{Proof of \Cref{main}}\label{proofmain}
	
	Let $u$ be a second Neumann eigenfunction. Suppose that
	$\nabla u(p_0)=0$ at some $p_0\in\Omega$. We choose a conformal
	map such that $\Psi(0)=p_0$, we set $W=J\nabla u$. The strict eigenvalue inequality
	\eqref{spectral-gap} therefore allows us to apply \Cref{hodge},
	which gives
	\begin{equation}\label{commonW}
		W\in\mathscr X,\qquad E_\mu(W)=0,\qquad
		E_\mu(Y)\ge0\quad \text{for all }Y\in\mathscr X.
	\end{equation}
	let $w$ be the disk coefficient of $W$, and set $F_W=e^\delta w$. 
	\Cref{lem:rotated-gradient} shows that $F_W$ satisfies
	\begin{equation}\label{rotated-second2}
		-4\pb(\kappa\pz F_W)=\mu\rho\kappa F_W
		\qquad\text{in }\D.
	\end{equation}
	Moreover, $F_W(0)=0$, and the $L^2$ mass and the Dirichlet energy of $u$ are given as
	\begin{subequations}\label{gradnorms2}
		\begin{align}
			\mu^2\|u\|_{L^2(\Omega)}^2
			&=4\int_\D\kappa|\pz F_W|^2 dz,	\label{rotated-derivative2}\\
			\|\nabla u\|_{L^2(\Omega)}^2
			&=\int_\D\rho\kappa|F_W|^2 dz.	\label{rotated-mass2}
		\end{align}
	\end{subequations}

	\subsection{The Euclidean plane}\label{plane}
	
	\begin{proof}[Proof of \Cref{main} in the Euclidean plane]
		In the Euclidean case, $\delta=0$, $\kappa=1$ and
		$\rho=|\Phi'|^2$. Since $W=J\nabla u=(-u_{x_2},u_{x_1})$,
		we identify its disk coefficient with a complex-valued function
		and choose
		\begin{equation}\label{plane-data}
			\Hh=\mathbb C,\qquad U=F_W=(-u_{x_2}+iu_{x_1})\circ\Psi,\qquad
			V(z)=\mu\rho(z)I_{\mathbb C}.
		\end{equation}
		Differentiating $-\Delta u=\mu u$ in the interior and then pulling
		the equations for the two components of $W$ back by $\Psi$ gives
		\begin{equation}\label{plane-equation}
			-\Delta U=\mu\rho U=VU\quad\text{in }\D,\qquad U(0)=0,
		\end{equation}
		where $U(0)=0$ follows from $\nabla u(p_0)=0$.
		
		The operator $V(z)$ is self-adjoint, nonnegative, bounded linear operator on the Hilbert space $\Hh$ for all $z\in\D$,
		and its trace $\operatorname{tr}V=\mu\rho\in L^r(\D)$ by \eqref{metricbounds}.
		For any $\xi\in H^1_0(\D;\Hh)$, since $\Psi$ induces an isomorphism between
		$H^1_0(\Omega)$ and $H^1_0(\mathbb D)$ via pullback,
		then $\Xi=\xi\circ\Psi^{-1}\in H^1_0(\Omega;\Hh)$, 
		Applying the Dirichlet Rayleigh inequality to the real and
		imaginary parts of $\Xi$, together with change of variables and conformal
		invariance of the Dirichlet integral, we obtain
		\begin{equation}\label{plane:gap}
			\int_\D\ip{\xi}{V\xi} dz
			=\mu\int_\Omega|\Xi|^2 dA
			\le\frac{\mu}{\lambda_D}\int_\Omega|\nabla\Xi|^2 dA
			=\frac{\mu}{\lambda_D}\int_\D\|\nabla\xi\|^2 dz .
		\end{equation}
		Since $\mu<\lambda_D$ holds by \Cref{lem:spectral-compare}, this verifies the assumption \eqref{small} with $\eta=\mu/\lambda_D<1$.

		To obtain the boundary data required by \Cref{strict},
		we decompose $U$ into an antiholomorphic part and an $H^1$ part.
		The derivative $p=\pz U$ is square integrable by \eqref{rotated-derivative2}
		and smooth in the interior, so we extend it by zero outside the disk
		and use the Cauchy transform
		\begin{equation}\label{cauchy}
			(\mathcal C p)(z)=\frac1\pi\int_\D\frac{p(t)}{\bar z-\bar t} dt,
			\qquad \pz\mathcal C p=p .
		\end{equation}
		The distributional derivative $\pz(1/(\pi\bar z))$ is the unit
		point mass at the origin, which gives the prescribed derivative.
		On the disk, the convolution kernel
		may be restricted to $|z-t|<2$, where its $L^1$ norm is $4$,
		so Young's inequality bounds $\mathcal Cp$ in $L^2$.
		For the other derivative, the Fourier symbols of $\pb$ and $\pz$
		are $\frac i2(\xi_1+i\xi_2)$ and $\frac i2(\xi_1-i\xi_2)$,
		whose ratio has modulus one away from the origin. Plancherel's
		theorem, first for smooth compactly supported functions and then
		by $L^2$ approximation, therefore yields
		\begin{equation}\label{cauchy-bound}
			\|\mathcal Cp\|_{L^2(\D;\Hh)}
			\le4\|p\|_{L^2(\D;\Hh)},\quad
			\|\pb\mathcal Cp\|_{L^2(\R^2;\Hh)}
			=\|p\|_{L^2(\D;\Hh)},\quad
			\|\mathcal Cp\|_{H^1(\D;\Hh)}
			\le C\|p\|_{L^2(\D;\Hh)}.
		\end{equation}
		Here the last $H^1$ bound follows by combining the two estimates
		with $|\nabla X|^2=2(|\pz X|^2+|\pb X|^2)$. Since $p$ is
		smooth in the interior, interior regularity for $\pz$ also makes
		$\mathcal Cp$ smooth in the interior, so we may normalize its value at the
		origin and set
		\begin{equation}\label{plane-splitting}
			\widehat U=\mathcal Cp-(\mathcal Cp)(0),\qquad U_{\rm a}=U-\widehat U.
		\end{equation}
		By construction, $\pz U_{\rm a}=0$ and $U_{\rm a}(0)=\widehat U(0)=0$, and \eqref{cauchy-bound} gives $\widehat U\in H^1(\D;\Hh)$.
		Identity \eqref{rotated-mass2} implies
		$\|\sqrt\rho U\|_{L^2(\D;\Hh)}=\|\nabla  u\|_{L^2(\Omega)}$, together with the weighted Sobolev estimate
		\eqref{weightspace}, gives
		\begin{equation}\label{plane-weighted}
			\|\sqrt\rho\,U_{\rm a}\|_{L^2(\D;\Hh)}
			\le\|\sqrt\rho\,U\|_{L^2(\D;\Hh)}
			+\|\sqrt\rho\,\widehat U\|_{L^2(\D;\Hh)}
			\le\|\nabla  u\|_{L^2(\Omega)}
			+C{\|\widehat U\|_{H^1(\D;\Hh)}} .
		\end{equation}
		Thus $S_VU_{\rm a}\in L^2(\D;\Hh)$, and taking
		$f=\widehat U_{\partial\D}\in H^{1/2}(\partial\D;\Hh)$ supplies
		the data for \eqref{base-problem}.
		Uniqueness of the solution of the equation identifies the constructed $U$ with its solution.

		Let $U_i$ be the solutions of \eqref{comparison-problems}, 
		let $Z_i$ be the vector fields on $\Omega$ such that
		$F_{Z_i}=U_i$.
		Then \eqref{field-local-identities-div-curl} and \eqref{solution-energy-data} give
		\begin{align}
			&	\|Z_i\|_{L^2(\Omega)}^2=\int_\D\rho|U_i|^2 dz = \mu^{-1}\|S_VU_i\|_{L^2(\D;\Hh)}^2 <\infty,\\
			&\|\diver Z_i\|_{L^2(\Omega)}^2+\|\curl Z_i\|_{L^2(\Omega)}^2
			=4\int_\D|\pz U_i|^2 dz<\infty,
		\end{align}
		Moreover, \Cref{boundary-multiplier} gives
		$U_i-\beta_iU\in H^1_0(\D;\Hh)$, with some $\beta_i$ real on the
		circle, so $F_{Z_i}-\beta_iF_W\in H^1_0(\D;\Hh)$ and \Cref{admissible} yields $Z_i\in\mathscr X$.
		We may therefore apply \Cref{diskenergy} to obtain
		\begin{equation}\label{plane-energy-identity}
			E_\mu(Z_i)=4\int_\D|\pz U_i|^2 dz-\mu\int_\D\rho|U_i|^2 dz
			=\mathcal A_V(U_i).
		\end{equation}
		For the original solution $U$, the mass identity \eqref{rotated-mass2}
		and the zero $E_\mu$ energy of $W$ give
		\begin{equation}\label{plane-mass}
			\|S_VU\|_{L^2(\D;\Hh)}^2
			=\mu\|\nabla u\|_{L^2(\Omega)}^2>0,\qquad
			\mathcal A_V(U)=E_\mu(W)=0 .
		\end{equation}
		Since $\|S_VU\|_{L^2(\D;\Hh)} >0$, \Cref{strict} gives
		$\mathcal A_V(U_1)+\mathcal A_V(U_2)-\mathcal A_V(U)<0$.
		The fields $Z_i$ belong to $\mathscr X$, so \Cref{hodge}
		and \eqref{plane-energy-identity} give
		\begin{equation}\label{plane-contradiction}
			\begin{aligned}
				0 \le E_\mu(Z_1)+E_\mu(Z_2) &=E_\mu(Z_1)+E_\mu(Z_2)-E_\mu(W) \\
				&=\mathcal A_V(U_1)+\mathcal A_V(U_2)-\mathcal A_V(U)<0 ,
			\end{aligned}
		\end{equation}
		which contradicts the assumed interior critical point and proves the theorem in $\mathbb R^2$.
	\end{proof}
	
	\subsection{The hyperbolic plane}\label{hyperbolic}
	
	In the hyperbolic case, we use the Poincar\'e coordinate $\zeta$. Let $\Omega_{\rm coord}=\zeta(\Omega)$
	and let $\Phi\colon\D\to\Omega_{\rm coord}$ be a Riemann map.
	Since $\Omega$ is bounded in the hyperbolic metric, 
	$\Omega_{\rm coord}$ is compactly contained in the unit disk,
	therefore $\sup_\D|\Phi(z)|\le R<1$ for some $R$. 
	In the notation of \eqref{conformaldata}, we have
	\begin{equation}\label{hyper-metric}
		\kappa(z)=e^{-2\delta(z)}=\frac{(1-|\Phi(z)|^2)^2}{4},\qquad
		\rho=\frac{4|\Phi'|^2}{(1-|\Phi|^2)^2}.
	\end{equation}

	We take $\mathcal H=\ell^2(\mathbb N_0)$ to be the Hilbert space of complex sequences,
	with its norm $\|(c_n)\|^2=\sum_n|c_n|^2$.
	Expanding $(1-|\Phi|^2)^{-2}$ by
	$(1-t)^{-2}=\sum_{n=0}^\infty(n+1)t^n$, we obtain
	\begin{equation}\label{hyper-column}
		\begin{aligned}
			b_n=2\sqrt{n+1}\,\Phi^n,\quad
			\mathbf{b}=(b_n)_{n\ge0},\quad \kappa^{-1}=\|\mathbf b\|^2=\sum_n  |b_n|^2.
		\end{aligned}
	\end{equation}
	Since $|\Phi|\le R<1$, for every $z\in\D$, we have
	\begin{align}\label{hyper-column-bound}
		\|\mathbf{b}(z)\|^2\le \sum_n 4(n+1)R^{2n}<\infty, \quad \|\nabla\mathbf{b}(z)\|^2 \le |\Phi'(z)|^2\sum_n 4n^2(n+1)R^{2n-2}<\infty.
	\end{align}
	Thus \eqref{metricbounds} implies that $\mathbf{b}(z)\in W^{1,2r}(\D;\mathcal H)\cap L^\infty(\D;\mathcal H)$.
	For each $z\in\D$, let $\mathbf b^*(z)=\overline{\mathbf b(z)}^T$ be a bounded linear functional
	from $\mathcal H$ to $\mathbb C$, defined by $\mathbf b^*(z) \mathbf  h =\sum_n\overline{b_n(z)}h_n$.
	Then it satisfies
	\begin{align}\label{hyper-row}
		|\mathbf b^*(z) \mathbf h|\le\|\mathbf{b}(z)\|\|\mathbf h\|= \kappa^{-1/2}(z)\|\mathbf h\|, \qquad  \pz\mathbf b^*(z)=0 .
	\end{align}
	
	If $X\in H^1(\D;\mathcal H)$, applying H\"older inequality and pointwise bound \eqref{hyper-row}, we have
	$\|\mathbf b^* X\|_{L^2(\D;\mathbb{C})}\le \|\mathbf b\|_{L^\infty(\D;\mathcal H)}\|X\|_{L^2(\D;\mathcal H)}$.
	The product rule yields the estimate
	$|\nabla(\mathbf b^* X)|\le \|\mathbf b\|\| \nabla X\| + \| \nabla \mathbf b\|\|X\| $,
	applying H\"older's inequality and the Sobolev embedding, notice that 
	$\mathbf{b}(z)\in W^{1,2r}(\D;\mathcal H)\cap L^\infty(\D;\mathcal H)$, we have
	\begin{equation}\label{rowbound}
		\begin{aligned}
			\|\nabla(\mathbf b^* X)\|_{L^2(\D )}
			&\le \|\mathbf b\|_{L^\infty(\D;\mathcal H)}
			\| \nabla X\|_{L^2(\D;\mathcal H)}
			+\|\nabla \mathbf b\|_{L^{2r}(\D;\mathcal H)}
			\|X\|_{L^{2r/(r-1)}(\D;\mathcal H)} \\
			&\le   C\|X\|_{H^1(\D;\mathcal H)}.
		\end{aligned}
	\end{equation}
	Therefore, we have $\mathbf b^* X \in H^1(\D)$.
	If $X\in H^1_0(\D;\mathcal H)$, then $\mathbf b^*X\in H^1_0(\D)$.

	Define $V(z)=\mu\rho(z)\kappa(z)\mathbf{b}(z) \mathbf{b}^*(z)$. Clearly, $V(z)$
	is nonnegative bounded self-adjoint operator on $\mathcal H$, 
	and for the standard basis $\mathbf e_n \in \mathcal H$, we have
	\begin{equation}\label{h:trace}
		\begin{aligned}
			\operatorname{tr}V(z)
			&=\sum_n\ip{\mathbf e_n}{V(z)\mathbf e_n}
			=\mu\rho(z) \kappa(z)\sum_n|b_n|^2=\mu\rho(z)\in L^r(\D).
		\end{aligned}
	\end{equation}
	This verify assumption \eqref{eq:trace-V}.
	Since $\|S_V \xi\|^2=\mu\rho \kappa|\mathbf b^*\xi|^2\le \mu\rho \|\xi\|^2$ by \eqref{hyper-row},
	analogous to the planar case, a change of variables and the conformal invariance of the Dirichlet integral yield
	\begin{equation}\label{hi:eq:potential}
		\begin{aligned}
			\int_\D\ip{\xi}{V\xi} dz
			=\mu\int_\D\rho\kappa|\mathbf b^*\xi|^2 dz
			\le \mu\int_\D\rho\|\xi\|^2 dz
			\le \frac{\mu}{\lambda_D}\int_\D\|\nabla\xi\|^2 dz \quad\text{ for }\xi\in H^1_0(\D;\mathcal H).
		\end{aligned}
	\end{equation}
	Since $\mu<\lambda_D$ holds by \Cref{lem:spectral-compare}, this verifies the assumption \eqref{small} with $\eta=\mu/\lambda_D<1$.

	\begin{lemma}\label{lifting}
		There exist an antiholomorphic $U_{\rm a}$ with
		$\sqrt\rho U_{\rm a}\in L^2(\D;\Hh)$ and a function $\widehat U\in H^1(\D;\Hh)$,
		both vanishing at $0$, such that $U=U_{\rm a}+\widehat U$ satisfies
		\begin{equation}\label{liftid}
			\pz U=\kappa\mathbf b \pz F_W,\qquad
			\mathbf b^*U=F_W,\qquad -\Delta U=VU .
		\end{equation}
		In addition,
		\begin{equation}\label{liftenergy}
			\mathcal A_V(U)=E_\mu(W)=0,\qquad
			\|S_VU\|_{L^2(\D;\Hh)}^2=\mu\|\nabla u\|_{L^2(\Omega)}^2>0 .
		\end{equation}
	\end{lemma}
	
	\begin{proof}
		Let $\mathbf p=\kappa\mathbf b \pz F_W$. Since $\kappa\|\mathbf b\|^2=1$, then \eqref{rotated-derivative2} gives
		\begin{equation}\label{hyper-source-norm}
			\sum_{n\ge0}\|p_n\|_{L^2(\D)}^2
			=\int_\D\kappa|\pz F_W|^2 dz=\tfrac14\mu^2\|u\|_{L^2(\Omega)}^2 .
		\end{equation}
		Thus $\mathbf p\in L^2(\D;\Hh)$, and is smooth in the interior.
		Now, we first seek $U$ satisfying
		\begin{equation}\label{ha:target}
			\partial_zU=\mathbf p,\qquad
			\mathbf b^* U=F_W,\qquad U(0)=0.
		\end{equation}
		We extend $\mathbf p$ by zero outside the disk and use the Cauchy transform
		\begin{align}
			(\mathcal C \mathbf p)(z)=\frac1\pi\int_\D\frac{\mathbf p(t)}{\bar z-\bar t} d t.
		\end{align}
		Analogous to the planar case, we have $\pz\mathcal C \mathbf p=\mathbf p $ and $\|\mathcal C\mathbf p\|_{H^1(\D;\Hh)}
		\le C\|\mathbf p\|_{L^2(\D;\Hh)}$.
		Therefore, we set
		$\widehat U=\mathcal C\mathbf p-(\mathcal C\mathbf p)(0)$, which satisfies
		$\pz\widehat U=\mathbf p$, $\widehat U(0)=0$ and $ \widehat U\in H^1(\D;\Hh)$.

		For the remaining constraint, define
		$\varepsilon=F_W-\mathbf b^*\widehat U$, then antiholomorphicity of
		$\mathbf b^*$ gives
		\begin{equation}\label{eq:varepsi-anti}
			\pz\varepsilon=(1-\kappa\|\mathbf b\|^2)\pz F_W=0,\qquad
			\varepsilon(0)=0.
		\end{equation}
		Thus $\varepsilon$ is antiholomorphic and vanishes at the origin.
		Moreover, the bound \eqref{rotated-mass2} for $F_W$, the weighted Sobolev inequality \eqref{weightspace},
		and \eqref{rowbound} give
		\begin{equation}\label{hyper-error-bound}
			\begin{aligned}
				\|\sqrt\rho {\varepsilon}\|_{L^2(\D)}
				&\le \|\sqrt\rho F_W\|_{L^2(\D)}+ \| \sqrt\rho \mathbf b^*\widehat U\|_{L^2(\D)} \\
				&\le C\|\nabla u\|_{L^2(\Omega)}+ C\|\mathbf b^*\widehat U\|_{H^1(\D)} \\
				&\le C\|\nabla u\|_{L^2(\Omega)}+C\|\widehat U\|_{H^1(\D;\Hh)} <\infty .
			\end{aligned}
		\end{equation}
		Notice  the $0$-th component of $\mathbf b$ is $b_0=2$ by \eqref{hyper-column},
		we set $U_{\rm a}=({\varepsilon}/2)\mathbf e_0$ and $U=\widehat U+U_{\rm a}$.
		Therefore, we have $U_{\rm a}(0)=0$, $\mathbf b^*U_{\rm a}={\varepsilon}$, $\mathbf b^*U
		=\mathbf b^*\widehat U+\mathbf b^*U_{\rm a}=F_W$,
		and \eqref{hyper-error-bound} gives $\sqrt{\rho} U_{\rm a}\in L^2(\D;\Hh)$.
		Moreover,  since ${\varepsilon}$ is antiholomorphic, we obtain
		$\partial_zU=\partial_z\widehat U=\mathbf p$.
		Applying $-4\pb$ to the  identity $\partial_zU=\mathbf p$  and using
		$\pb\mathbf b=0$, together with \eqref{rotated-second2}, we have
		\begin{equation*}
			-\Delta U=-4\pb\partial_zU =-4\mathbf b\,\pb(\kappa\pz F_W)
			=\mu\rho\kappa\mathbf b F_W=VU .
		\end{equation*}
		This proves the disk equation. Finally, 
		$\|\mathbf b\|^2=\kappa^{-1}$ and  $V=\mu\rho\kappa\mathbf{b}\mathbf{b}^*$ give
		\begin{equation*}
			\|\pz U\|^2=\kappa|\pz F_W|^2,\qquad
			\|S_VU\|^2=\mu\rho \kappa|\mathbf b^*U|^2=\mu\rho\kappa|F_W|^2 .
		\end{equation*}
		Therefore, integrating these identities, by use \Cref{diskenergy} and \Cref{hodge}, together with $W\in \mathscr{X}$, we have
		\begin{equation}
			\begin{aligned}
				\mathcal A_V(U)&=4\int_\D\|\partial_zU\|^2 dz -\int_\D\|S_VU\|^2 dz\\
				&=4\int_\D\kappa|\partial_zF_W|^2 dz -\mu\int_\D\rho\kappa|F_W|^2 dz\\
				&=E_\mu(W)=0.
			\end{aligned}
		\end{equation}
		and \eqref{rotated-mass2} gives
		\begin{align}
			\|S_VU\|_{L^2(\D;\Hh)}^2=\mu\int_\D\rho\kappa|F_W|^2 dz=\mu\|\nabla u\|_{L^2(\Omega)}^2>0,
		\end{align}
		where the last inequality is positive because $u$ is a nonzero Neumann eigenfunction.
	\end{proof}
	
	\begin{proof}[Proof of \Cref{main} in the hyperbolic plane]
		The function $U=U_{\rm a}+\widehat U$ constructed in \Cref{lifting} satisfies $\sqrt\rho U_{\rm a}\in L^2(\D;\Hh)$,
		then \eqref{hyper-row} gives
		\[
		\|S_VU_{\rm a}\|_{L^2(\D;\Hh)}^2
		=\mu\int_\D\rho\kappa|\mathbf b^*U_{\rm a}|^2 dz
		\le\mu\int_\D\rho\|U_{\rm a}\|^2 dz<\infty.
		\]
		Taking $f=\widehat U|_{\partial\D}\in H^{1/2}(\partial\D;\Hh)$,
		uniqueness of equation \eqref{base-problem} identifies this $U$ with its
		solution.
		Let $U_i$ be the solutions of \eqref{comparison-problems}, 
		let $Z_i$ be the vector fields on $\Omega$ such that
		$F_{Z_i}=\mathbf b^* U_i$.
		Since $\pz\mathbf b^*=0$, then $\pz F_{Z_i}= \mathbf b^* \pz U_i $.
		Therefore, \eqref{hyper-row} and the identity $\|S_V \xi\|^2=\mu\rho \kappa|\mathbf b^*\xi|^2$ gives
		\begin{equation}\label{project}
			\kappa|\pz F_{Z_i}|^2\le\|\pz U_i\|^2,\qquad
			\|S_VU_i\|^2=\mu\rho\kappa|F_{Z_i}|^2.
		\end{equation}
		Therefore, \eqref{field-local-identities-div-curl}, \eqref{project}  and \eqref{solution-energy-data} give
		\begin{equation}
			\begin{aligned}
				\|Z_i\|_{L^2(\Omega)}^2
				&=\int_\D\rho\kappa|F_{Z_i}|^2 dz
				=\mu^{-1}\|S_VU_i\|_{L^2(\D;\Hh)}^2<\infty,\\
				\|\diver Z_i\|_{L^2(\Omega)}^2
				+\|\curl Z_i\|_{L^2(\Omega)}^2
				&=4\int_\D\kappa|\pz F_{Z_i}|^2 dz
				\le4\|\pz U_i\|_{L^2(\D;\Hh) }^2<\infty.
			\end{aligned}
		\end{equation}
		Moreover, \Cref{boundary-multiplier} gives
		$U_i-\beta_iU\in H^1_0(\D;\Hh)$, with some $\beta_i$ real on the
		circle,
		since $F_W=\mathbf b^*U$ and $F_{Z_i}=\mathbf b^* U_i$, we have
		\begin{equation}\label{hyper-boundary}
			F_{Z_i}-\beta_i F_W=\mathbf b^*(U_i-\beta_i U)\in H^1_0(\D).
		\end{equation}
		\Cref{admissible} therefore gives $Z_i\in\mathscr X$.
		Using \Cref{diskenergy} and \eqref{project}, we obtain
		\begin{equation}\label{hyper-projected-energy}
			\begin{aligned}
				\mathcal A_V(U_i)-E_\mu(Z_i)
				&=\left(4\int_\D\|\pz U_i\|^2 dz
				-\int_\D\|S_VU_i\|^2 dz\right)
				-\left(4\int_\D\kappa|\pz F_{Z_i}|^2 dz
				-\mu\int_\D\rho\kappa|F_{Z_i}|^2 dz\right)\\
				&=4\int_\D\bigl(\|\pz U_i\|^2
				-\kappa|\mathbf b^*\pz U_i|^2\bigr) dz\ge0.
			\end{aligned}
		\end{equation}
        By \eqref{liftenergy}, $\mathcal A_V(U)=0$ and $\|S_VU\|_{\mathcal L}>0$.
		Combining \Cref{strict} with \eqref{hyper-projected-energy} and
		\eqref{commonW}, we obtain
		\begin{equation}\label{hyper-contradiction}
			0\le E_\mu(Z_1)+E_\mu(Z_2)
			\le\mathcal A_V(U_1)+\mathcal A_V(U_2)<\mathcal A_V(U)=0 .
		\end{equation}
		This contradiction proves the hyperbolic assertion.
	\end{proof}
	
	\subsection{The unit sphere}\label{sphere}
	
	In the unit sphere case, we use the  stereographic
	projection $\zeta$ from the pole opposite the center of the containing hemisphere on the sphere. 
	Since $\Omega$ is compactly contained in an open hemisphere,
	$\Omega_{\rm coord}$ is compactly contained in the unit disk,
	therefore $\sup_\D|\Phi(z)|\le R<1$ for some $R$. 
	In the notation of \eqref{conformaldata}, we have
	\begin{equation}\label{sphere-metric}
		m(\zeta)=\frac4{(1+|\zeta|^2)^2},\qquad
		\kappa(z)=e^{-2\delta(z)}=\frac{(1+|\Phi|^2)^2}{4}
		,\qquad
		\rho=\frac{4|\Phi'|^2}{(1+|\Phi|^2)^2}.
	\end{equation}
	By using \eqref{rotated-derivative2} and \eqref{rotated-mass2},
	together with that $\kappa$ is bounded above and below, we have
	$\pz F_W,\sqrt\rho\,F_W\in L^2(\D)$. Thus the Cauchy transform construction
	used in the planar case gives
	$\widehat F_W=\mathcal C(\pz F_W)-(\mathcal C(\pz F_W))(0)$,
	the estimate \eqref{cauchy-bound}, together with
\eqref{rotated-mass2}, the positive lower bound for $\kappa$
and the weighted Sobolev estimate \eqref{weightspace}, yields the decomposition
	\begin{equation}\label{sphsplit}
		\begin{aligned}
			&F_W=F_{W,\rm a}+\widehat F_W,\qquad
			F_{W,\rm a}(0)=\widehat F_W(0)=0,\\
			&\pz F_{W,\rm a}=0,\qquad \sqrt\rho\,F_{W,\rm a}\in L^2(\D),\qquad
			\widehat F_W\in H^1(\D).
		\end{aligned}
	\end{equation}
	
	Set
	\begin{equation}\label{sphere-coefficients}
		F=e^{-\delta}w=\kappa F_W .
	\end{equation}
	We now derive the equation for $F$. A direct computation for $\kappa$ in \eqref{sphere-metric} gives
	\begin{equation}\label{sphere-log-kappa}
		\pz\pb\log\kappa
		=\frac{2|\Phi'|^2}{(1+|\Phi|^2)^2}
		=\frac\rho2 .
	\end{equation}
	Dividing \eqref{rotated-second2} by $\kappa$ gives
	$	-4\bigl(\pz\pb F_W+(\pb\log\kappa)\pz F_W\bigr)=\mu\rho F_W$, together with
	$\kappa^{-1}\pb(\kappa F_W)=\pb F_W+(\pb\log\kappa)F_W$, differentiate this identity in $z$ and use
	\eqref{sphere-log-kappa} to obtain
	\begin{equation}
		\begin{aligned}
	-4\pz\bigl(\kappa^{-1}\pb(\kappa F_W)\bigr)=\mu\rho F_W-4(\pz\pb\log\kappa)F_W	=(\mu-2)\rho F_W.
	\end{aligned}
	\end{equation}
	Thus $F$ satisfies
	\begin{equation}\label{sphshift}
		-4\pz(\kappa^{-1}\pb F)
		=\tau\rho\kappa^{-1}F,
		\qquad \tau=\mu-2.
	\end{equation}

	Take $\Hh=\mathbb C^3$, the identity $\kappa=(1+|\Phi|^2)^2/4$ in \eqref{sphere-metric}
	can be decomposed into the sum of three terms $1/4$, $|\Phi|^2/2$ and $|\Phi|^4/4$. 
	Based on this identity, we set
	\begin{equation}\label{sphere-column}
		\mathbf b=\tfrac12(1,\sqrt2\,\Phi,\Phi^2)^T,\quad
		\|\mathbf b\|^2=\kappa,\quad
		\Pi=\kappa^{-1}\overline{\mathbf b}\mathbf b^T, \quad V=\tau\rho\Pi,
	\end{equation}
	where $\mathbf b^T h=\sum b_j h_j$ is a bounded linear function from $\Hh$ to $\mathbb{C}$,
	and $\Pi$ is a projection on $\Hh$.
	Clearly, $\Pi^2=\Pi=\Pi^*$, and $\operatorname{tr} V=\tau\rho\in L^r(\D)$,
	this verify the assumption \eqref{eq:trace-V}.
	Since $\|S_V\xi\|^2=\tau\rho \|\Pi^2\xi\|^2$,
	analogous to the previous case, a change of variables and the conformal invariance of the Dirichlet integral yield
	\begin{equation}\label{sphere-gap}
		\int_\D\|S_V\xi\|^2 dz=\tau\int_\D\rho\|\Pi\xi\|^2 dz\le \tau\int_\D\rho\|\xi\|^2 dz
		\le\frac{\tau}{\lambda_D}\int_\D|\nabla\xi|^2 dz.
	\end{equation}
	Since $\mu<\lambda_D$ holds by \Cref{lem:spectral-compare}, together with 
	$\tau=\mu-2>0$, this verifies the assumption \eqref{small} with
	$\eta=\tau/\lambda_D<1$.
	
	We shall also need the fact, analogous to \eqref{rowbound},
	that multiplication by the coefficients preserves $H^1$.
	In fact, \eqref{metricbounds} and the boundedness of $\kappa$ imply that
	$\mathbf{b},\kappa$, $\kappa^{-1}$ and $\Pi$ all are belong to $W^{1,2r}\cap L^\infty$.
	The product rule implies that boundedness of multiplication on $H^1$, and approximation by compactly supported
	functions shows that it also preserves $H^1_0$.

	\begin{lemma}\label{lem:spherical-lifting}
		Set $U_{\rm a}=\overline{\mathbf b} F_{W,\rm a}$. There is a
		$\widehat U\in H^1(\D;\Hh)$ such that $U=U_{\rm a}+\widehat U$
		satisfies
		\begin{equation}\label{sphlift}
			\mathbf b^T U=F,\qquad
			\pb U={\kappa^{-1}} \overline{\mathbf b} \pb F,\qquad
			-\Delta U=V U,\qquad U(0)=0 .
		\end{equation}
		Moreover $\pz U=\pz\widehat U\in L^2(\D;\Hh)$ and
		$(I-\Pi)\pb U=0$ in $\D$.
	\end{lemma}
	\begin{proof}
		First observe that the first two identities in \eqref{sphlift} imply the third equation.
		Indeed,
		since $\pz\overline{\mathbf b}=0$, \eqref{sphshift} gives
		\begin{equation}\label{sphere-equation}
			-\Delta U=-4\overline{\mathbf b} \pz( \kappa^{-1}\pb F)
			=\tau\rho\kappa^{-1} \overline{\mathbf b}F=\tau\rho\Pi U =VU.
		\end{equation}
		It therefore suffices to satisfy the first two identities
		in \eqref{sphlift}. We seek $U=U_{\rm a}+\widehat U$, where
		$U_{\rm a}=\overline{\mathbf b}F_{W,\rm a}$ is a fixed function,
		hence the required identity for $\partial_{\bar z}U$ then reduces to
		\[
		\partial_{\bar z}\widehat U
		=\kappa^{-1}\overline{\mathbf b}\,\partial_{\bar z}F
		-\partial_{\bar z}U_{\mathrm{a}}=:\mathbf p.
		\]
		Now we compute $\mathbf p$.
		Differentiating $\kappa=\mathbf b^T\overline{\mathbf b}$ gives $\pb\kappa=\mathbf b^T\pb\overline{\mathbf b}$
		and $\kappa^{-1}\overline{\mathbf b}\pb\kappa=\Pi\pb\overline{\mathbf b}$.
		Using $F=\kappa(F_{W,\rm a}+\widehat F_W)$ and
		$U_{\rm a}=\overline{\mathbf b}F_{W,\rm a}$, we have
		\begin{equation}\label{sphere-source}
		\begin{aligned}
		\mathbf p&:=\kappa^{-1}\overline{\mathbf b} \pb F-\pb U_{\rm a}=\kappa^{-1}\overline{\mathbf b} \pb(\kappa F_{W,\rm a})
		  +\kappa^{-1}\overline{\mathbf b} \pb(\kappa\widehat F_W)-\pb(\overline{\mathbf b}F_{W,\rm a})\\
		&=\bigl(\kappa^{-1}\overline{\mathbf b} \pb\kappa  -\pb\overline{\mathbf b}\bigr)F_{W,\rm a}
		  +\kappa^{-1}\overline{\mathbf b}   \pb(\kappa\widehat F_W)
		=-(I-\Pi)\pb\overline{\mathbf b} F_{W,\rm a} +\kappa^{-1}\overline{\mathbf b}\pb(\kappa\widehat F_W).
		\end{aligned}
		\end{equation}
		Since $|(I-\Pi)\pb\overline{\mathbf b}|\le\|\mathbf b'\|\le C|\Phi'|\le C\sqrt\rho$,
		and multiplication by the coefficients $\kappa$  preserves $H^1$,
		together with the decomposition \eqref{sphsplit}, we have
		\begin{equation*}
			\|\mathbf p\|_{L^2(\D;\Hh)}\le C\|\sqrt\rho\,F_{W,\rm a}\|_{L^2(\D)}
			+C\|\widehat F_W\|_{H^1(\D)}<\infty .
		\end{equation*}

		Using the conjugate Cauchy transform, set
		$U_0=\overline{\mathcal C(\mathbf{\bar p})}
		-(\overline{\mathcal C(\mathbf {\bar p})})(0)$. Similarly to \eqref{cauchy-bound}, we have $U_0\in H^1(\D;\Hh)$,
		$\pb U_0=\mathbf p$ and $U_0(0)=0$.
		For the remaining row error, define $\varepsilon=F-\mathbf b^T(U_{\rm a}+U_0)$,
		the decomposition \eqref{sphsplit} gives $\varepsilon=\kappa\widehat F_W-\mathbf b^TU_0$,
		then $\varepsilon\in H^1(\D)$ by the multiplier bound, and we have

		\begin{equation*}
			\pb\varepsilon
			=\pb F-\mathbf b^T(\pb U_{\rm a}+ \mathbf p)
			=(1-\kappa^{-1}\|\mathbf b\|^2)\pb F=0,\qquad \varepsilon(0)=0.
		\end{equation*}
		Thus $\varepsilon$ is holomorphic and vanishes at the origin,
		so adding $(\varepsilon/b_0)\mathbf e_0$ to $U_0$ corrects the
		row without changing its $\pb$ derivative, where $b_0=1/2$.
		Therefore, we set 
		\begin{align}
			\widehat U=U_0+\frac{\varepsilon}{b_0}\mathbf e_0,\quad  U=U_{\rm a}+\widehat U.
		\end{align}
		Thus $\widehat U\in H^1(\D;\Hh)$ proves the first two identities
		of \eqref{sphlift} and hence the third equation, with $U_{\rm a}(0)=\widehat U(0)=U(0)=0$. 
		Lastly, $U_{\rm a}$ is antiholomorphic, so we also have
		\begin{equation}\label{sphere-orthogonal}
			(I-\Pi)\pb U=0,\qquad \pz U=\pz\widehat U\in L^2(\D;\Hh).
		\end{equation}
	\end{proof}

Let $U_{\rm a}, \widehat{U}$ and $U=U_{\rm a}+\widehat U$ are the function constructed in \Cref{lem:spherical-lifting}.	
	Since $U_{\rm a}=\overline{\mathbf b}F_{W,\rm a}$ and $\|\mathbf b\|^2=\kappa$, \eqref{sphsplit} yields
	\begin{align}
		\|S_VU_{\rm a}\|_{L^2(\D;\Hh)}^2
	=\tau\int_\D\rho\kappa|F_{W,\rm a}|^2 dz<\infty
	\end{align}
The function $\widehat U$ belongs to $H^1(\D;\Hh)$, so its boundary value
$f=\widehat U|_{\partial\D}$ lies in $H^{1/2}(\partial\D;\Hh)$,
uniqueness of equation \eqref{base-problem}
identifies $U=U_{\rm a}+\widehat U$ with its solution. Let $U_i$ be the solutions of
\eqref{comparison-problems}.

For each $Z=U,U_1,U_2$,
we define its projected component and its orthogonal component by
\begin{equation}\label{sphere-decomposition}
	h_Z=\kappa^{-1}\mathbf b^TZ,\qquad
	T_Z=(I-\Pi)Z,\qquad
	Z=\overline{\mathbf b} h_Z+T_Z.
\end{equation}
Since $\Pi Z=\overline{\mathbf b} h_Z$, this is indeed an orthogonal decomposition.
We now determine the regularity of these components, write the three solutions in the form
\begin{equation}\label{sphere-trial-split}
	\begin{aligned}
		U&=\overline{\mathbf b} F_{W,\rm a}+\widehat U,\\
		U_i&=\overline{\mathbf b} (\overline{\Gamma_i}F_{W,\rm a})
		+(v_{i,\rm r}+q_{v_i}).
	\end{aligned}
\end{equation}
In either case, $Z$ can be written $Z=\overline{\mathbf b}D+Z_{\rm r}$,
where $D$ is antiholomorphic and $Z_{\rm r}\in H^1(\D;\Hh)$.
Moreover, $\sqrt\rho D\in L^2(\D)$.
Indeed, if $Z=U$, $\sqrt\rho D\in L^2(\D)$ is \eqref{sphsplit}, 
and if $Z=U_i$, the pole in
$D=\overline{\Gamma_i}F_{W,\rm a}$ is removable by
\eqref{removable-products}, since $\kappa$ is bounded below, the estimate \eqref{singular-mass}
then gives $\sqrt\rho\,D\in L^2(\D)$. Thus $\sqrt\rho D\in L^2(\D)$ holds.
Because $\Pi\overline{\mathbf b}=\overline{\mathbf b}$ and
$\mathbf b^T\overline{\mathbf b}=\kappa$,
substituting $Z=\overline{\mathbf b}D+Z_{\mathrm{r}}$ into \eqref{sphere-decomposition} yields
\begin{equation}\label{sphere-regularity}
	\begin{aligned}
		T_Z&=(I-\Pi)(\overline{\mathbf b}D+Z_{\rm r})=(I-\Pi)Z_{\rm r}\in H^1(\D;\Hh),\\
		h_Z&=\kappa^{-1}\mathbf b^T\overline{\mathbf b}\,D
		+\kappa^{-1}\mathbf b^TZ_{\rm r}=D+\kappa^{-1}\mathbf b^TZ_{\rm r},
		\qquad \pz h_Z,\ \sqrt\rho\,h_Z\in L^2(\D),
	\end{aligned}
\end{equation}
where the integrability follows from the $H^1$ multiplier bound for $I-\Pi$ and $\kappa^{-1}\mathbf b^T$, together with the weighted Sobolev estimate \eqref{weightspace}.

Let $Y[Z]$ be the tangent field whose disk coefficient is $e^{-\delta}h_Z$,
then \eqref{field-local-identities-div-curl} and \eqref{sphere-regularity} gives
\[
\|Y[Z]\|_2^2=\int_\D\rho\kappa|h_Z|^2 dz<\infty,\qquad
\|\diver Y[Z]\|_2^2+\|\curl Y[Z]\|_2^2
=4\int_\D\kappa|\pz h_Z|^2 dz<\infty.
\]
Clearly, for $Z=U$, the identity $\mathbf b^TU=F=\kappa F_W$ in \eqref{sphlift} gives $h_U=F_W$, and therefore $Y[U]=W$.
Next we verify the boundary condition for $Y[U_i]$.
\Cref{boundary-multiplier} gives
$U_i-\beta_iU\in H^1_0(\D;\Hh)$, with some $\beta_i$ real on the
circle,
apply the $H^1$ multiplier bound $\kappa^{-1}\mathbf b^T$ to
$U_i-\beta_iU$, together with $h_U=F_W$, we have
\begin{align}\label{eq:xxow}
	h_{U_i}-\beta_iF_W
=\kappa^{-1}\mathbf b^T(U_i-\beta_iU)
\in H^1_0(\D),
\end{align}
where the multiplier $\beta_i$ is real on $\partial\D$, so
\Cref{admissible} yields $Y[U_i]\in\mathscr X$.
For $Z=U$, we already have $Y[U]=W\in\mathscr X$ by
\eqref{commonW}. \Cref{diskenergy} therefore applies to
all three fields and gives
\begin{align}\label{eq:e-mu-tangent}
	E_\mu(Y[Z])
=4\int_\D\kappa|\pz h_Z|^2 dz
-\mu\int_\D\rho\kappa|h_Z|^2 dz,
\qquad Z=U,U_1,U_2.
\end{align}
Here $E_\mu(Y[Z])$ are finite because \eqref{sphere-regularity}.

\begin{lemma}\label{lem:spherical-comparison}
		For $Z=U,U_1,U_2$, the difference
		\begin{equation}\label{sphnull}
			\begin{aligned}
				\mathcal B(Z)
				&:=\mathcal A_V(Z)-E_\mu(Y[Z])-4\|(I-\Pi)\pb Z\|_2^2\\
				&=4\bigl(\|\pz T_Z\|_2^2-\|\pb T_Z\|_2^2\bigr)
				-8\operatorname{Re}\int_\D
				\pb\bigl(\overline{h_Z}\,(\mathbf b')^T T_Z\bigr) dz
			\end{aligned}
		\end{equation}
		is finite and satisfies
		\begin{equation}\label{sphere-cancellation}
			\mathcal B(U_1)+\mathcal B(U_2)=\mathcal B(U).
		\end{equation}
		Moreover $Y[U_i]\in\mathscr X$, and
		\begin{equation}\label{sphbridge}
			\mathcal A_V(U_1)+\mathcal A_V(U_2)-\mathcal A_V(U)
			=\sum_{i=1}^2 E_\mu(Y[U_i])
			+4\sum_{i=1}^2\|(I-\Pi)\pb U_i\|_2^2\ge0 .
		\end{equation}
	\end{lemma}
	\begin{proof}
		Write $h=h_Z$, $T=T_Z$, and
		$\mathbf t=(I-\Pi)\pb\overline{\mathbf b}$. The definition of $\Pi$ and $\mathbf b^T\overline{\mathbf b}=\kappa$
		gives $\mathbf b^TT=\mathbf b^T(I-\Pi)Z=0$. Differentiating
		$\mathbf b^TT=0$ and $Z=\overline{\mathbf b}h+T$,
		together with $\pb\mathbf b=0$, we have
		\begin{equation}\label{sphere-proj-derivatives}
			\begin{aligned}
				&\mathbf b^T\pz T=\pz(\mathbf b^TT)-(\mathbf b')^TT=-(\mathbf b')^TT,\qquad\qquad\qquad\qquad\quad
				\mathbf b^T\pb T=\pb(\mathbf b^TT)=0,\\
				&\pz Z=\overline{\mathbf b} \pz h+\pz T,\qquad\qquad\ 
				(I-\Pi)\pb Z =(I-\Pi)(\pb\overline{\mathbf b})h +(I-\Pi)\pb T=\mathbf t h+\pb T .
			\end{aligned}
		\end{equation}
		Here the last equality uses $\Pi\pb T=0$ follows from
		$\mathbf b^T\pb T=0$. 
		Expanding the squared norms in the last two identities of \eqref{sphere-proj-derivatives},
		and using $\mathbf t^*=(\mathbf b')^T(I-\Pi)$, we obtain
		\begin{equation}\label{sphere-expansion}
			\begin{aligned}
				&|\pz Z|^2
				=|\overline{\mathbf b} \pz h|^2+|\pz T|^2
				  +2\operatorname{Re}\bigl(\overline{\pz h} \mathbf b^T\pz T\bigr)
				=\kappa|\pz h|^2+|\pz T|^2
				-2\operatorname{Re}\bigl(\overline{\pz h} (\mathbf b')^TT\bigr),\\
				&|(I-\Pi)\pb Z|^2
				=|\mathbf t|^2|h|^2+|\pb T|^2
				  +2\operatorname{Re}\bigl(\bar h \mathbf t^*\pb T\bigr)
				=|\mathbf t|^2|h|^2+|\pb T|^2
				+2\operatorname{Re}\bigl(\bar h (\mathbf b')^T\pb T\bigr).
			\end{aligned}
		\end{equation}
		and by using $\kappa=\mathbf b^T\overline{\mathbf b}$, the squared norm of $\mathbf t$ is gives by
		\begin{equation}\label{sphere-curvature}
			\begin{aligned}
			|\mathbf t|^2
			=|\pb\overline{\mathbf b}|^2
			  -|\Pi\pb\overline{\mathbf b}|^2
			&=\|\mathbf b'\|^2
			  -\kappa^{-1}|\mathbf b^T\pb\overline{\mathbf b}|^2\\
			&=\pz\pb\kappa
			  -\kappa^{-1}(\pz\kappa)(\pb\kappa)
			 =\kappa\pz\pb\log\kappa
			 =\tfrac12\kappa\rho .
			\end{aligned}
		\end{equation}
		Since $\Pi Z=\overline{\mathbf b}h$ and $|\mathbf b|^2=\kappa$, we have
		\[
		\begin{aligned}
			\mathcal A_V(Z)
			&=4\|\partial_zZ\|_{L^2(\D;\Hh)}^2-\|S_VZ\|_{L^2(\D;\Hh)}^2
			=4\int_{\mathbb D}|\partial_zZ|^2 dz
			-\tau\int_{\mathbb D}\rho|\Pi Z|^2 dz\\
			&=4\int_{\mathbb D}|\partial_zZ|^2 dz
			-\tau\int_{\mathbb D}\rho\kappa|h|^2 dz.
		\end{aligned}
		\]
		Together with \eqref{eq:e-mu-tangent} and \eqref{sphere-expansion}, we obtain
		\[
		\begin{aligned}
		\mathcal B(Z)
		&=4\int_\D\bigl(|\pz Z|^2-\kappa|\pz h|^2
		              -|(I-\Pi)\pb Z|^2\bigr) dz
		 +(\mu-\tau)\int_\D\rho\kappa|h|^2 dz\\
		&=4\bigl(\|\pz T\|_2^2-\|\pb T\|_2^2\bigr)
		 +\int_\D\bigl((\mu-\tau)\rho\kappa-4|\mathbf t|^2\bigr)|h|^2 dz
		-8\operatorname{Re}\int_\D
		\bigl(\overline{\pz h}(\mathbf b')^TT
		+\bar h(\mathbf b')^T\pb T\bigr) dz.
		\end{aligned}
		\]	
		Here $\mu-\tau=2$, so the coefficient of $|h|^2$ vanishes
		by \eqref{sphere-curvature}. Moreover, $\overline{\pz h}(\mathbf b')^TT
		+\bar h(\mathbf b')^T\pb T
		=\pb\bigl(\bar h(\mathbf b')^TT\bigr)$
		because $\pb\mathbf b'=0$. This proves \eqref{sphnull}. Clearly $\mathcal B(Z)$ are finite.

		We next compare the boundary values. By \eqref{tracecorrection},
		multiplication by $I-\Pi$ and $\kappa^{-1}\mathbf b^T$ gives
		\begin{equation}\label{sphere-differences}
			\begin{aligned}
			T_{U_i}-\beta_i T_U
			&=(I-\Pi)(U_i-\beta_iU)\in H^1_0(\D;\Hh),\\
			h_{U_i}-\beta_i h_U
			&=\kappa^{-1}\mathbf b^T(U_i-\beta_iU)\in H^1_0(\D).
			\end{aligned}
		\end{equation}
		Consider first the squared-derivative difference in \eqref{sphnull},
		for $X=X_1+iX_2\in H^1(\D;\mathbb C)$, expanding
		$2\pz X$ in real coordinates gives
		\begin{equation}\label{jacobian-identity}
			4|\pz X|^2=|\nabla X|^2+2\det D(X_1,X_2),\qquad
			|\pz X|^2-|\pb X|^2=\det D(X_1,X_2).
		\end{equation}
		For smooth $X$, the Stokes' formula gives
		\begin{equation}\label{jacobian-boundary}
			\begin{aligned}
				\int_\D (|\pz X|^2-|\pb X|^2) dz
				=\tfrac12\int_{\partial\D}(X_1\mathrm dX_2-X_2\mathrm dX_1).
			\end{aligned}
		\end{equation}
		The left-hand side is continuous under $H^1$ convergence, while
		the boundary expression is the continuous pairing of an $H^{1/2}$
		value with a tangential derivative in $H^{-1/2}$, so the identity
		extends to $H^1$ and depends only on $X_{\partial\D}$.
		For a real smooth multiplier $\gamma$, the product rule gives
		\begin{equation*}
			(\gamma u_1)\mathrm d(\gamma u_2)-(\gamma u_2)\mathrm d(\gamma u_1)
			=\gamma^2(u_1\mathrm du_2-u_2\mathrm du_1).
		\end{equation*}
		The identity extends to $H^{1/2}$ boundary data by approximation.
		Writing $u=(T_U)|_{\partial\D}$ componentwise, we have
		$(T_{U_i})|_{\partial\D}=\gamma_i u$ by
		\eqref{sphere-differences}, hence the boundary integrands satisfy
		\[
		\sum_{i=1}^2
		\bigl((\gamma_i u_1)\mathrm d(\gamma_i u_2)
		      -(\gamma_i u_2)\mathrm d(\gamma_i u_1)\bigr)
		=(\gamma_1^2+\gamma_2^2)(u_1\mathrm du_2-u_2\mathrm du_1)
		=u_1\mathrm du_2-u_2\mathrm du_1.
		\]
		Applying \eqref{jacobian-boundary} to each of the three
		components of $T$ yields
		\begin{equation}\label{jacobian-sum}
			\sum_{i=1}^2\big(\|\pz T_{U_i}\|_2^2-\|\pb T_{U_i}\|_2^2\big)
			=\|\pz T_U\|_2^2-\|\pb T_U\|_2^2 .
		\end{equation}
		
		For the remaining term set
		\begin{equation}\label{mixed-form}
			\begin{aligned}
			I(h,T) =\int_\D
			\bigl(\overline{\pz h} (\mathbf b')^TT
			+\bar h (\mathbf b')^T\pb T\bigr) dz
			=\int_\D\pb\bigl(\bar h(\mathbf b')^TT\bigr) dz .
			\end{aligned}
		\end{equation}
		For $h,T$ with the regularity in \eqref{sphere-regularity}, by using H\"older inequality and Sobolev embedding, the
		bound $|\mathbf b'|\le C|\Phi'|\le C\sqrt\rho$ yields
		\begin{equation}\label{mixed-bound}
			\begin{aligned}
				&\left|\int_\D\overline{\pz h}(\mathbf b')^TT dz\right|
				\le C\|\pz h\|_{L^2(\D)}\|\rho\|_{L^r(\D)}^{1/2}\|T\|_{L^{2r/(r-1)}(\D;\Hh)}
				\le C\|\pz h\|_{L^2(\D)}\|\rho\|_{L^r(\D)}^{1/2}\|T\|_{H^1(\D;\Hh)},\\
				&\left|\int_\D\bar h(\mathbf b')^T\pb T dz\right|\le C\|\sqrt\rho h\|_{L^2(\D)}\|\pb T\|_{L^2(\D;\Hh)} .
			\end{aligned}
		\end{equation}
		Thus both integrals are finite, including when $h$ contains an antiholomorphic part
		with only weighted $L^2$ control.
		If one factor is smooth and compactly supported, the last integral
		in \eqref{mixed-form} vanishes. For $h_j\to h$ in $H^1_0(\D)$,
		\eqref{weightspace} and \eqref{mixed-bound} give
	\begin{align}
			|I(h_j-h,T)|\le C
		\bigl(\|\pz(h_j-h)\|_{L^2(\D)}
		     +\|\sqrt\rho (h_j-h)\|_{L^2(\D)}\bigr)\|T\|_{H^1(\D;\Hh)}
		\longrightarrow0.
	\end{align}
		Likewise, if $T_j\to T$ in $H^1_0(\D;\Hh)$, then
	\begin{align}
			|I(h,T_j-T)|\le C
		\bigl(\|\pz h\|_{L^2(\D)}+\|\sqrt\rho h\|_{L^2(\D)}\bigr)
		\|T_j-T\|_{H^1(\D;\Hh)}\longrightarrow0.
	\end{align}
		Consequently $I(h,T)=0$ if either $h$ or $T$ belongs to $H^1_0$
		and the other satisfies \eqref{sphere-regularity}.

		Since $\beta_i$ and its gradient are bounded and $\pz\beta_i$
		is supported away from the boundary, the products
		$\beta_i h_U$ and $\beta_i T_U$ satisfy the respective bounds
		for $h_U$ and $T_U$ in \eqref{sphere-regularity}.
		Using \eqref{sphere-differences}, $h_{U_i}$ can be decomposed into $\beta_i h_U$ plus an $H^1_0$ part, and $T_{U_i}$ can be decomposed into $\beta_i T_U$ plus an $H^1_0$ part. By the linearity of $I$ and the fact that $I$ vanishes whenever one of the factors is in $H^1_0$, we have
		\begin{equation}\label{mixed-cancellation}
			\begin{aligned}
				\sum_{i=1}^2 I(h_{U_i},T_{U_i})-I(h_U,T_U)
				&=\sum_{i=1}^2 I(\beta_i h_U,\beta_i T_U)-I(h_U,T_U)\\
				&=\int_\D\pb\!\left[
				\bar h_U(\mathbf b')^T
				\bigl((|\beta_1|^2+|\beta_2|^2-1)T_U\bigr)
				\right] dz\\
				&=I\bigl(h_U,(|\beta_1|^2+|\beta_2|^2-1)T_U\bigr)=0 ,
			\end{aligned}
		\end{equation}
		where the last expression vanishes because
		$(|\beta_1|^2+|\beta_2|^2-1)T_U\in H^1_0$. Together with \eqref{jacobian-sum} and \eqref{mixed-cancellation},
		we have
	\begin{equation}
			\begin{aligned}
			\mathcal B(U_1)+\mathcal B(U_2)-\mathcal B(U)
			=&4\left[
			\sum_{i=1}^2\bigl(\|\partial_zT_{U_i}\|_2^2	-\|\partial_{\bar z}T_{U_i}\|_2^2\bigr)-\bigl(\|\partial_zT_U\|_2^2
			-\|\partial_{\bar z}T_U\|_2^2\bigr)
			\right]\\
			&-8\operatorname{Re}\left[\sum_{i=1}^2 I(h_{U_i},T_{U_i})-I(h_U,T_U)\right]=0.
		\end{aligned}
	\end{equation}
		This proves \eqref{sphere-cancellation}.
		
		We have already verified that $Y[U_i]\in\mathscr X$,
		From \eqref{sphnull} and \eqref{sphere-cancellation}, we obtain
		\begin{equation}
			\begin{aligned}
		\mathcal A_V(U_1)+\mathcal A_V(U_2)-\mathcal A_V(U)
		&=E_\mu(Y[U_1])+E_\mu(Y[U_2])-E_\mu(Y[U])\\
		&+4\sum_{i=1}^2\|(I-\Pi)\pb U_i\|_2^2-4\|(I-\Pi)\pb U\|_2^2 .
		\end{aligned}
		\end{equation}
		Now $Y[U]=W$ and $E_\mu(W)=0$ by \eqref{commonW}, while
		$(I-\Pi)\pb U=0$ by \eqref{sphere-orthogonal}. Finally,
		$E_\mu(Y[U_i])\ge0$ by \Cref{hodge}. This proves
		\eqref{sphbridge}.
		
	\end{proof}

\begin{proof}[Proof of \Cref{main} on the unit sphere]
Let $\mu>2$ and put $\tau=\mu-2$.  We use the Hilbert space $\Hh=\mathbb C^3$ and the
potential $V=\tau\rho\Pi$ from \eqref{sphere-column}. The estimate
\eqref{sphere-gap} verifies the form bound required in
\Cref{strict}. By \Cref{lem:spherical-lifting}, the disk
solution $U=U_{\rm a}+\widehat U$ satisfies $-\Delta U=VU$ and $U(0)=0$.
The weighted bound for $F_{W,\rm a}$ in \eqref{sphsplit} gives
$S_VU_{\rm a}\in L^2(\D;\Hh)$. Taking
$f=\widehat U_{\partial\D}$ therefore identifies $U$ with the
solution of \eqref{base-problem}. Let $U_1,U_2$ be the corresponding
solutions of \eqref{comparison-problems}.

For each $i$, decompose $U_i$ into its projected component $h_{U_i}$ and its orthogonal component,
and take the field $Y_i=Y[U_i]$ defined by the disk
coefficient $e^{-\delta}h_{U_i}$ in \eqref{sphere-decomposition}.
The bounds in \eqref{sphere-regularity} and the local identities
\eqref{field-local-identities-div-curl} show that
$Y_i,\diver Y_i,\curl Y_i\in L^2(\Omega)$. Moreover,
\eqref{tracecorrection} and the multiplication bound above give
\[
 h_{U_i}-\beta_iF_W
   =\kappa^{-1}\mathbf b^T(U_i-\beta_iU)\in H^1_0(\D).
\]
Since $\beta_i$ is real on the circle, \Cref{admissible} yields
$Y_i\in\mathscr X$. Meanwhile, for the original solution, we have $Y[U]=W\in\mathscr X$.

We now compare the three disk energies. Rearranging \eqref{sphnull}
for $Z=U,U_1,U_2$ and taking the paired difference gives
\[
\begin{aligned}
 Q_V
 &=\sum_{i=1}^2\left(E_\mu(Y_i)
      +4\|(I-\Pi)\pb U_i\|_2^2+\mathcal B(U_i)\right)\\
 &\quad-\left(E_\mu(W)
      +4\|(I-\Pi)\pb U\|_2^2+\mathcal B(U)\right)\\
 &=\sum_{i=1}^2 E_\mu(Y_i)
      +4\sum_{i=1}^2\|(I-\Pi)\pb U_i\|_2^2\ge0.
\end{aligned}
\]
Here the terms $\mathcal B(U_i)$ cancel with $\mathcal B(U)$ by
\eqref{sphere-cancellation}, $E_\mu(W)=0$ by \eqref{commonW},
and $(I-\Pi)\pb U=0$ by \eqref{sphere-orthogonal}. 
\Cref{hodge} gives
$E_\mu(Y_i)\ge0$, and the other terms are squared norms.
Thus $Q_V\ge 0$, this is precisely \eqref{sphbridge}.
Finally, $\Pi U=\overline{\mathbf b}F_W$ and
$\|\mathbf b\|^2=\kappa$, so \eqref{rotated-mass2} gives
\begin{equation}\label{sphere-mass}
 \|S_VU\|_{L^2(\D;\Hh)}^2
 =\tau\int_\D\rho\kappa|F_W|^2 dz
 =\tau\|\nabla u\|_{L^2(\Omega)}^2>0.
\end{equation}
The last inequality holds because $u$ is a nonzero eigenfunction and
$\tau>0$. Since $U(0)=0$, \Cref{strict} now yields $Q_V<0$,
contradicting the nonnegativity above. This proves the spherical assertion.
\end{proof}

\noindent\textbf{AI disclosure:} We used ChatGPT (OpenAI) to help us calculate the strict two-multiplier inequality on the disk and derive the orthogonal decomposition of the lifted fields in the spherical case. We also used it to check the grammar throughout the manuscript.
The authors have reviewed all AI-assisted content and take full
responsibility for the mathematical results and the final manuscript.

\end{document}